\documentclass[11pt]{article}

\usepackage{amscd}      
\usepackage{amssymb}     
\usepackage[intlimits]{amsmath}     
\usepackage{xypic}      
\LaTeXdiagrams          
\usepackage[all,v2]{xy}      
\xyoption{2cell}
\UseAllTwocells
\xyoption{frame}
\CompileMatrices

\usepackage{theorem}

\newtheorem{prop}{Proposition}[section]  
\newtheorem{lem}[prop]{Lemma}

\newtheorem{cor}[prop]{Corollary}
\newtheorem{them}[prop]{Theorem}

\theorembodyfont{\upshape}

\newtheorem{defn}[prop]{Definition}

\newtheorem{numrmk}[prop]{Remark}

\newtheorem{numex}[prop]{Example}

\newenvironment{pf}{\begin{trivlist}\item[]{\sc Proof.}}%
            {\nolinebreak $\Box$ \end{trivlist}}

\newcommand{\noprint}[1]{}

\renewcommand{\tilde}{\widetilde}

\newcommand{\toto}{\rightrightarrows}

\newcommand{\upst}{^{\ast}}

\newcommand{\com}{^{\scriptscriptstyle\bullet}}
\newcommand{\lcom}{_{\scriptscriptstyle\bullet}}
\newcommand{\upcom}{^{\scriptscriptstyle\bullet}}

\newcommand{\XX}{{\mathfrak X}}

\newcommand{\Gg}{{\mathfrak g}}

\renewcommand{\O}{{\cal O}}

\newcommand{\lL}{{\cal L}}

\newcommand{\del}{\partial}

\newcommand{\pr}{\mathop{\rm pr}\nolimits}

\renewcommand{\Im}{\mathop{\rm Im}}

\newcommand{\Ad}{\mathop{\rm Ad}\nolimits}
\newcommand{\rank}{\mathop{\rm rank}\nolimits}

\newcommand{\ol}{\overline}

\newcommand{\ra}{\rangle}
\newcommand{\la}{\langle}   
\newcommand{\ldiag}[1]%
       {\makebox[0cm]{${\scriptstyle#1}\downarrow\phantom{\scriptstyle#1}$}}
\newcommand{\ldiagup}[1]%
       {\makebox[0cm]{${\scriptstyle#1}\uparrow\phantom{\scriptstyle#1}$}}
\newcommand{\rdiag}[1]%
       {\makebox[0cm]{$\phantom{\scriptstyle#1}\downarrow{\scriptstyle#1}$}}
\newcommand{\sediagr}[1]%
       {\makebox[0cm]{$\phantom{\scriptstyle#1}\searrow{\scriptstyle#1}$}}
\newcommand{\nediagr}[1]%
       {\makebox[0cm]{$\phantom{\scriptstyle#1}\nearrow{\scriptstyle#1}$}}
\newcommand{\rdiagup}[1]%
       {\makebox[0cm]{$\phantom{\scriptstyle#1}\uparrow{\scriptstyle#1}$}}
\newcommand{\swdiag}[1]%
       {\makebox[0cm]{$\phantom{\scriptstyle#1}\swarrow{\scriptstyle#1}$}}
\newcommand{\sediag}[1]%
       {\makebox[0cm]{${\scriptstyle#1}\searrow\phantom{\scriptstyle#1}$}}
\newcommand{\nediag}[1]%
       {\makebox[0cm]{${\scriptstyle#1}\nearrow\phantom{\scriptstyle#1}$}}

\newcommand{\longiso}{\stackrel{\textstyle\sim}{\longrightarrow}}

\newcommand{\hol}{\mbox{Hol}}

\newcommand{\doublearrowstack}[2]%
                      {{{{\scriptstyle#1}\atop{\textstyle\longrightarrow}}\atop{{\textstyle\longrightarrow}\atop{\scriptstyle#2}}}}
\newcommand{\rightleftarrowstack}[2]%
                      {{{{\scriptstyle#1}\atop{\textstyle\longrightarrow}}\atop{{\textstyle\longleftarrow}\atop{\scriptstyle#2}}}}
\newcommand{\leftrightarrowstack}[2]%
                      {{{{\scriptstyle#1}\atop{\textstyle\longleftarrow}}\atop{{\textstyle\longrightarrow}\atop{\scriptstyle#2}}}}

\newcommand{\overtoparrow}%
{\makebox[0cm]{\beginpicture
\setcoordinatesystem units <.8cm,.4cm> point at 0 0
\setplotarea x from -3 to 3, y from 0 to 1
\setquadratic
\plot -3 0 0 1 3 0 /
\put{\vector(3,-1){0}}[Bl] at 3 0
\endpicture}}

\newcommand{\underbottomarrow}%
{\makebox[0cm]{\beginpicture
\setcoordinatesystem units <.8cm,.4cm> point at 0 0
\setplotarea x from -3 to 3, y from 0 to 1
\setquadratic
\plot -3 1 0 0 3 1 /
\put{\vector(3,1){0}}[Bl] at 3 1
\endpicture}}

\newcommand{\ses}[5]%
{0\longrightarrow#1\stackrel{#2}{ \longrightarrow}#3\stackrel{#4}{
\longrightarrow}#5\longrightarrow0}

\newcommand{\dt}[6]%
{#1\stackrel{#2}{longrightarrow}#3 \stackrel{#4}{\longrightarrow}#5
\stackrel{#6}{\longrightarrow} #1[1]}  
 
\newcommand{\cat}[1]%
{(\mbox{\rm #1})}

\newcommand{\gm}{\Gamma }            
\newcommand{\lon }{\longrightarrow } 
\newcommand{\backl}{\mathbin{\vrule width1.5ex height.4pt\vrule height1.5ex}}
\newcommand{\per}{\backl}
\newcommand{\be}{\begin{eqnarray*}}
\newcommand{\ee}{\end{eqnarray*}}
\newcommand{\smalcirc}{\mbox{\tiny{$\circ $}}}

\newcommand{\half}{\frac{1}{2}}
\newcommand{\frakg}{\mathfrak g}
\newcommand{\frakh}{\mathfrak h}

\newcommand{\bimodule}{bimodule}
\newcommand{\tJ}{\tilde{J}} 
\newcommand{\btheta}{\bar{\theta}}

\newcommand{\complex}{{\mathbb{C}}}
\newcommand{\mmm}{J(x)}

\let\Vec=\overrightarrow
\let\ceV=\overleftarrow
\let\Hat=\widehat

\title{Momentum  Maps and Morita Equivalence}
\author{Ping Xu
         \thanks{ Research partially supported by NSF
       grants DMS00-72171 and DMS03-06665. }\\
        Department of Mathematics\\
         Pennsylvania State University \\
         University Park, PA 16802, USA\\
{\sf email: ping@math.psu.edu }
}     
 
\date{}
\begin{document}
\sloppy
\maketitle

\centerline{{Dedicated to Alan Weinstein on the occasion of his 60th birthday}}

\begin{abstract}
We introduce quasi-symplectic groupoids and explain their relation with momentum
map theories. This approach enables us to unify into a single framework various
momentum map theories, including ordinary Hamiltonian $G$-spaces, Lu's momentum
maps of Poisson group actions, and the group-valued momentum maps of
Alekseev--Malkin--Meinrenken. More precisely, we carry out the following
program: 

(1) We define and study properties of quasi-symplectic groupoids.

(2) We study the momentum map theory defined by a quasi-symplectic groupoid $\gm
\toto P$. In particular, we study the reduction theory and prove that
$J^{-1}(\O)/\gm $ is a symplectic manifold for any Hamiltonian $\gm$-space
$(X\stackrel{J}{\to}P, \omega_{X})$ (even though $\omega_{X} \in \Omega^2 (X)$
may be degenerate), where $\O\subset P$ is a groupoid orbit. More generally, we
prove that the intertwiner space $(X_1 \times_{P} \overline{X_{2}})/\gm$ between
two Hamiltonian $\gm$-spaces $X_1$ and $X_2$ is a symplectic manifold
(whenever it is a smooth manifold).

(3) We study Morita equivalence of quasi-symplectic groupoids. In particular, we
prove that Morita equivalent quasi-symplectic groupoids give rise to equivalent
momentum map theories. Moreover the intertwiner space
 $(X_1 \times_{P} \overline{X_{2}})/\gm$  depends only on
  the Morita equivalence class. As a result, we
recover various well-known results concerning equivalence of momentum maps
including the Alekseev--Ginzburg--Weinstein linearization theorem and the
Alekseev--Malkin--Meinrenken equivalence theorem between quasi-Hamiltonian
spaces and Hamiltonian loop group spaces.
\end{abstract}

\tableofcontents

\section{Introduction}

``Momentum'' usually refers to quantities whose conservation under the time
evolution of a physical system is related to some symmetry of the system.
Noether \cite{Noether}, in the course of developing ideas of Einstein and Klein
in general relativity theory, found a very general equivalence between
symmetries and conservation laws in field theory; this is now known as Noether's
theorem. Focusing on the relation between symmetries and conserved quantities,
the study of momentum maps has received much attention in the last three
decades, continuing to the present day with the formulation of new notions of
symmetry. In geometric terms, a phase space with a symmetry group consists of a
symplectic (or Poisson) manifold $P$ and an Hamiltonian action of a Lie group
$G$. By the latter, we mean a symplectic (or Poisson) action of $G$ on $P$
together with an equivariant map $J: P \to \frakg^*$ such that for each $X \in
\frakg$, the one-parameter group of transformations of $P$ generated by $X$ is
the flow of the Hamiltonian vector field with Hamiltonian $\la J(x), X\ra \in
C^{\infty}(P)$. The map $J$ is called the momentum (or moment) map of the
Hamiltonian action. One very important aspect of the momentum map theory is the
study of Marsden--Weinstein (or symplectic) reduction, which is the simultaneous
use of symmetries and conserved quantities to reduce the dimension of a
Hamiltonian system.

With the advance of physics and mathematics, new notions of symmetry and
momentum have appeared. For instance, a Poisson group symmetry is the classical
limit of a ``quantum group symmetry'' in quantum group theory \cite{dr:quantum}.
Lu's momentum map theory \cite{Lu:1991} for Poisson Lie group actions is a
theory adapted from the usual Hamiltonian theory which incorporates the Poisson
structure on the symmetry group $G$. Computations of the symplectic structures
on moduli spaces of flat connections on surfaces have led to another notion of
Hamiltonian symmetry known as quasi-Hamiltonian symmetry. In this new theory,
the $2$-form $\omega$ on the phase space is neither closed nor non-degenerate,
but these ``defects'' are compensated for by the presence of an auxiliary
structure on the group. This is the starting point of the theory of
quasi-Hamiltonian $G$-spaces with group-valued momentum maps of
Alekseev--Malkin--Meinrenken (AMM) \cite{AMM}. All these momentum map theories
share many similarities, but involve different techniques and proofs. It is also
known that some of these momentum theories are equivalent to one another. For
instance, for compact groups, the AMM group-valued momentum map theory is
equivalent to the Hamiltonian momentum map theory of loop groups of
Meinrenken--Woodward \cite{MeinW, MeinW1, MeinW2}, and for compact
Bruhat--Poisson groups, Lu's momentum map theory is equivalent to the usual
Hamiltonian momentum map theory \cite{A}. However, these results are fragmentary
and their geometric significance remains unclear. It is therefore natural to
investigate the relations between these theories, and to seek a uniform
framework, which is an open question raised by Weinstein \cite{Weinstein:02}.
A unified approach would seek to develop a single momentum map theory which
reduces to the theories already established under special circumstances. While
necessarily generalizing the problem, this would allow a direct comparison of
the features of the various momentum maps in a more intrinsic manner. The
importance of such a single momentum map theory is not merely to give another
interpretation of these existing momentum map theories, but rather to explore
the intrinsic ingredients of these theories so that techniques in one theory can
be applied to another. This is particularly important in the study of
group-valued momentum map theory where there are still many open problems,
including the quantization problem which we believe will be the main application
of our approach \cite{LX}.
 
The approach taken in this paper involves extending the notion of symmetry from
actions of groups to actions of groupoids. This was motivated by the work of
Mikami--Weinstein \cite{MW} who showed that the usual Hamiltonian momentum map
is in fact equivalent to the symplectic action of the symplectic groupoid $T^*G
\toto \frakg^*$, which integrates the Lie-Poisson structure on $\frakg^*$.
Similarly, in \cite{WeinsteinX:yang}, 
 Weinstein and the author proved 
that  the 
momentum map theory of Lu for an Hamiltonian  Poisson group $G$-space
is equivalent to the symplectic action of
 the symplectic groupoid $G\times G^*\toto G^*$
integrating the dual Poisson group $G^*$ \cite{LuW:1989}. By a symplectic action
of a symplectic groupoid $\gm \toto P$ on a symplectic manifold $X$, we mean a
map $J: X\to P$ equipped with a $\gm$-action $\gm \times_{P}X \to X$ which is
compatible with the symplectic structures \cite{MW}. In this case $X$ is called
an Hamiltonian-$\gm$ space.

There is strong evidence that the AMM group-valued momentum map is closely
related to the transformation groupoid $G\times G\toto G$. Here $G$ acts on
itself by conjugation. However, $G\times G\toto G$ is no longer a symplectic
groupoid since the closed $3$-form, i.e., the Cartan form $\Omega$ on $G$, must
now play a role. In fact, one can show that the standard AMM $2$-form $\omega
\in \Omega^2(G \times G)$ together with $\Omega \in \Omega^3(G)$ gives a
$3$-cocycle of the total de Rham complex of the groupoid and defines a
nontrivial class in the equivariant cohomology $H^3_{G}(G)$ \cite{BXZ}.

This example suggests that one must enrich the notion of a symplectic groupoid
in order to include such ``twisted'' symplectic structures on the groupoids.
Thus we arrive at quasi-symplectic groupoids, the main subject of the present
paper. A quasi-symplectic groupoid is a Lie groupoid ${\gm}\toto{P}$ equipped
with a $2$-form $\omega \in \Omega^2(\gm)$ and a $3$-form $\Omega \in
\Omega^3(P)$ such that $\omega + \Omega$ is a $3$-cocycle of the de Rham complex
of the groupoid, where $\omega$ must satisfy a weak non-degeneracy condition.
When $\omega$ is honestly non-degenerate, this is the so-called twisted
symplectic groupoid studied by Cattaneo and the author \cite{CX} as the global
object integrating a twisted Poisson structure of Severa--Weinstein \cite{SW}.
In particular, when $\Omega$ vanishes, it reduces to an ordinary symplectic
groupoid.

It turns out that much of the theory of Hamiltonian $\gm$-spaces of a symplectic
groupoid $\gm$ can be generalized to the present context of quasi-symplectic
groupoids. In particular, one can perform reduction and prove that
$J^{-1}(\O)/\gm $ is a symplectic manifold (even though $\omega_{X} \in
\Omega^2(X)$ may be degenerate), where $\O \subset P$ is an orbit of the
groupoid. More generally, one can introduce the classical intertwiner space
$(X_1 \times_{P} \overline{X_{2}})/\gm$ between two Hamiltonian $\gm$-spaces
$X_1$ and $X_2$, generalizing the same notion studied by Guillemin--Sternberg
\cite{GS} for the ordinary Hamiltonian $G$-spaces. One shows that this is a
symplectic manifold (whenever it is a smooth manifold).

As for symplectic groupoids, one can also introduce Morita equivalence for
quasi-symplectic groupoids. In particular, we prove the following main result.
(i) Morita equivalent quasi-symplectic groupoids give rise to equivalent
momentum map theories in the sense that there is an equivalence of categories
between their Hamiltonian $\gm$-spaces; (ii) the symplectic manifold $(X_1
\times_{P} \overline{X_{2}})/\gm$  depends only on
 the Morita equivalence class
of $\gm$. As a result, we recover various well-known results concerning
equivalence of momentum maps including the Alekseev--Ginzburg--Weinstein
linearization theorem and Alekseev--Malkin--Meinrenken equivalence theorem for
group-valued momentum maps. They are essentially due to the Morita equivalence
between the Lu--Weinstein symplectic groupoid $G \times G^* \toto G^*$ and the
standard cotangent symplectic groupoid $T^*G\toto \frakg^*$, where $G$ is a
compact simple Lie group equipped with the Bruhat--Poisson group structure and
the Morita equivalence is between the symplectic groupoid $(LG\times L\frakg
\toto {L\frakg}, \omega_{LG \times L\frakg})$ and the AMM quasi-symplectic
groupoid $(G\times G\toto G, \omega+\Omega)$.

Another main motivation of the present work is the quantization problem. It is
natural to study the geometric quantization of the symplectic reduced space
$J^{-1}(\O)/\gm$ or more generally the symplectic intertwiner space $(X_1
\times_{P} \overline{X_{2}})/\gm$, and prove the Guillemin--Sternberg conjecture
that ``$[Q, R]=0$'' for Hamiltonian $\gm$-spaces. As an application, our uniform
framework naturally leads to the following construction of prequantizations. A
prequantization of the quasi-symplectic groupoid $(\gm \toto P, \omega + 
\Omega)$ is a gerbe over the stack corresponding to the groupoid $\gm\toto P$,
while a prequantization of an Hamiltonian $\gm$-space is a line bundle $L$ on
which the gerbe acts. A prequantization of the symplectic intertwiner space
$(X_1 \times_{P} \overline{X_{2}})/\gm$ can be constructed using these data. For
symplectic groupoids, such a prequantization was studied in \cite{Xu:intert}.
Details of this construction for quasi-symplectic groupoids  appear
elsewhere \cite{LX}. Note that in the usual Hamiltonian case, since the
symplectic $2$-form defines a zero class in the third cohomology group of the
groupoid $T^{*}G\toto \frakg^*$, which in this case is the equivariant
cohomology $H^3_G(\frakg^*)$, gerbes do not enter explicitly. However, for a
general quasi-symplectic groupoid (for instance the AMM quasi-symplectic
groupoid), since the $3$-cocycle $\omega + \Omega$ may define a nontrivial
class, gerbes are inevitable in the construction.

Recently, Zung proved the convexity theorem for Hamiltonian $\gm$-spaces
 of proper quasi-symplectic groupoids, which encompasses many classical convexity theorems
in the literature \cite{Zung}. 
Finally we note that recently Bursztyn--Crainic--Weinstein--Zhu showed that
infinitesimally quasi-symplectic groupoids (which are called twisted
presymplectic groupoids in \cite{BCWZ}) correspond to twisted Dirac structures.
They also studied the infinitesimal version of our Hamiltonian $\gm$-spaces. We
refer to \cite{BCWZ} for details.

{\bf Acknowledgments.}
The author
 would like to thank several institutions for their hospitality while work on
this project was being done: RIMS/Kyoto University, Ecole Polytechnique, Erwin
Schr\"odinger Institute, and University of Geneva. He also wishes to thank many
people for useful discussions and comments, including Anton Alekseev, Philip
Boalch, Henrique Bursztyn, Eckhard Meinrenken, Jim Stasheff, and Alan Weinstein.
Some results of the paper were presented in Poisson 2002, Lisbon. He
 would like to thank the organizers for
 inviting him and giving him the chance to present
the work. Special thanks go to Laurent-Gengoux family, who generously provided a
fantastic working environment during the S\'eminaire Itin\'erant G\'eom\'etrie
et Physique I, Normandie 2003, where this work was completed.

\section{Quasi-symplectic groupoids}

In this section, we introduce quasi-symplectic groupoids and discuss
their basic properties.

\subsection{Pre-quasi-symplectic  groupoids}

A simple and compact way to define a pre-quasi-symplectic  groupoid
is to use the de-Rham double complex of a Lie groupoid.
First, let us recall its definition  below.

Let $\Gamma\toto \gm_{0}$ be a Lie groupoid with source and
target maps $s, t: \gm \to \gm_{0}$. Define for all
 $p\geq0$
$$\Gamma_{p}= \underbrace{\Gamma\times_{\gm_{0}}\ldots\times_{\gm_{0}}\Gamma}_{\text{$p$
times}}\,,$$
i.e., $\Gamma_{p}$ is the manifold of composable 
sequences of $p$ arrows in the groupoid $\Gamma\toto \gm_{0}$.
We have $p+1$ canonical maps $\gm_{p}\to \gm_{p-1}$
giving rise to a              diagram
\begin{equation}\label{sim.ma}
\xymatrix{
\ldots \gm_{2}
\ar[r]\ar@<1ex>[r]\ar@<-1ex>[r] & \gm_{1}\ar@<-.5ex>[r]\ar@<.5ex>[r]
&\gm_{0}\,.}
\end{equation}

 
In fact,  $\gm \lcom$ is a simplicial manifold.
Consider  the double complex $\Omega\upcom(\Gamma\lcom)$:
\begin{equation}
\label{eq:DeRham}
\xymatrix{
\cdots&\cdots&\cdots&\\
\Omega^1(\gm_{0})\ar[u]^d\ar[r]^\partial &\Omega^1(\gm_{1})\ar[u]^d\ar[r]^\partial
&\Omega^1(\gm_{2})\ar[u]^d\ar[r]^\partial&\cdots\\
\Omega^0(\gm_{0})\ar[u]^d\ar[r]^\partial&\Omega^0(\gm_{1})\ar[u]^d\ar[r]^\partial
&\Omega^0(\gm_{2})\ar[u]^d\ar[r]^\partial&\cdots
}
\end{equation}                           
Its boundary maps are $d:
\Omega^{k}( \gm_{p} ) \to \Omega^{k+1}( \gm_{p} )$, the usual exterior
derivative of differentiable forms and $\partial
:\Omega^{k}( \gm_{p} ) \to \Omega^{k}( \gm_{p+1} )$,  the
alternating sum of the pull-back maps of (\ref{sim.ma}).
We denote the total differential by $\delta = (-1)^pd+\del$.
The  cohomology groups of the total complex $\Omega\upcom(\Gamma\lcom)$
$$H_{DR}^k(\Gamma\lcom)= H^k\big(\Omega\upcom(\Gamma\lcom)\big)$$
are called the {\em de~Rham cohomology }groups of $\Gamma\toto \gm_{0}$.
We now introduce the   following

\begin{defn}
A pre-quasi-symplectic  groupoid is  a Lie groupoid  ${\gm}\toto{P}$  
equipped with a two-form $\omega \in \Omega^2 (\gm )$ and
a three-form $\Omega\in \Omega^3 (P)$ such that
\begin{equation}
d \Omega =0, \ \   d \omega=\partial \Omega, \ \ \mbox{and } \  \partial \omega =0.
\end{equation}
In other words, $\omega+\Omega $ is  a 3-cocycle of the  total 
de-Rham complex of the groupoid  $\gm\toto P$.
\end{defn}     

\begin{numrmk}
It is simple to see that the last condition  $\partial \omega =0$
is equivalent to that the graph of the multiplication
$\Lambda \subset \gm \times  \gm \times \overline{\gm}$
is isotropic. In this case, $\omega$ is said  to be
multiplicative.
\end{numrmk}

By $A\to P$ we denote the Lie algebroid of $\gm\toto P$, where the
anchor map is denoted by $a:A\to TP$.
For any $\xi \in \gm (A)$, by $\Vec{\xi} $ and $\ceV{\xi}$
we denote its corresponding  right-and left-invariant vector fields on
$\gm$ respectively.
The following properties can be easily verified (see also \cite{CX}).

\begin{prop}
\label{prop:2.2}
Let $(\gm \toto P, \omega+\Omega)$ be  a pre-quasi-symplectic 
groupoid.
\begin{enumerate}
\item $\epsilon^* \omega= 0$,
 where $\epsilon :P\to \gm$ is the unit map;
\item $i^* \omega= -\omega$, where $i:\gm \to \gm$ is the  groupoid
inversion;
\item  for any $\xi, \eta \in \gm (A)$, 
$$\omega (\Vec{\xi} , \Vec{\eta})=
-\omega (\ceV{\xi} ,\ceV{\eta}), \ \ \ \ \omega (\Vec{\xi} ,\ceV{\eta})=0;$$
\item    for any $\xi, \eta \in \gm (A)$, $\omega (\Vec{\xi} , \Vec{\eta})$
is  a right invariant function on $\gm$,
 and $\omega (\ceV{\xi} ,\ceV{\eta})$ is a 
left invariant function on $\gm$.                               
\end{enumerate}
\end{prop} 
\begin{pf}
Let $\Lambda =\{(x, y, z)|z=xy, \ (x, y)\in \gm_2\}\subset \gm \times 
\gm \times  \overline{\gm}$ be the graph of groupoid multiplication.
Thus $\Lambda$ is isotropic with respect to $(\omega, 
\omega, -\omega)$.

(1).  For any $\delta'_{m}, \ \delta''_{m}\in T_m P$,
since $(\delta'_{m} , \delta'_{m} , \delta'_{m}),\
(\delta''_{m} , \delta''_{m}, \delta''_{m}) \in T\Lambda $,
it follows that $\omega (\delta'_{m} , \delta''_{m} )=0$.

(2). $\forall x\in \gm$ 
and $\forall \delta'_{x}, \ \delta''_{x}\in T_{x}\gm$,
 it is clear that  
$ (\delta'_{x} , i_{*} \delta'_{x} , s_* \delta'_{x} ),
 \ (\delta''_{x} , i_{*} \delta''_{x} , s_* \delta''_{x} )
\in T\Lambda $. Thus using (1),
 we have $$\omega (\delta'_{x}, \delta''_{x})+\omega (i_{*} \delta'_{x},
i_{*} \delta''_{x} )=0, $$ and therefore (2) follows.

(3).   Since $i_*\Vec{\xi}=-\ceV{\xi}$ and
$i_*\Vec{\eta}=-\ceV{\eta}$, from (2) it follows that
$\omega (\Vec{\xi} , \Vec{\eta})=
-\omega (\ceV{\xi} ,\ceV{\eta})$.
Now for any $x\in \gm$,  since both vectors
$( \Vec{\xi}(x), 0_{t (x)}, \Vec{\xi}(x))$  and
$(0_x , \ceV{\eta}(t (x) ), \ceV{\eta}(x ))$ are tangent to
$\Lambda $, we thus have $\omega (\Vec{\xi}(x), \ceV{\eta}(x ))=0$.
 
(4). It is simple to see that,
 for any $\xi, \eta \in \gm (A)$ and any composable pair
$(x, y)\in \gm_2$,
 $( \Vec{\xi}(x), 0_{y}, \Vec{\xi}(xy))$, 
$( \Vec{\eta}(x), 0_{y}, \Vec{\eta }(xy)) \in T\Lambda$.
Thus
$$\omega (\Vec{\xi}(x),  \Vec{\eta}(x))-\omega (\Vec{\xi}(xy) ,
\Vec{\eta }(xy) )=0. $$
Hence $\omega (\Vec{\xi} , \Vec{\eta})$
is  a right invariant function on $\gm$.  Similarly,
  one proves that
$\omega (\ceV{\xi} ,\ceV{\eta})$ is a left invariant function
 on $\gm$.            
\end{pf}

We next investigate the kernel of $\omega$ along the
unit  space $P$. For any $m\in P$,  there are two ways
to identify  elements of $A_{m}$ as tangent vectors of
$\gm$, namely vectors tangent to the $t$-fiber $\xi \to \Vec{\xi}(m)$,
or  to the $s$-fiber $\xi \to \ceV{\xi}(m)$. Write

\begin{equation}
 \Vec{A}|_{m}=\{\Vec{\xi}(m)|\forall \xi \in A_{m}\},\ \ \ \ 
\mbox{and }\ceV{A}|_{m}=\{\ceV{\xi}(m)|\forall \xi \in A_{m}\}.
\end{equation} 

Thus we have the following
decomposition of the tangent space:
\begin{equation}
T_{m}\gm =\Vec{A}|_{m}\oplus T_{m}P =\ceV{A}|_{m}\oplus T_{m}P , \ \ 
\forall m\in P.
\end{equation}

\begin{cor} 
\label{cor:2.3}
Under the same hypothesis as in Proposition \ref{prop:2.2}, we have, 
for any $m\in P$,
\begin{enumerate}
\item $\ker \omega_m =(\ker \omega_m \cap  \Vec{A}|_{m} ) \oplus
(\ker \omega_m \cap T_{m}P )$;
\item if $\Vec{\xi}(m) \in \ker \omega_m$, then 
$a (\xi )\in \ker \omega_m$;  and
\item for any $\xi \in A_{m}$, $\Vec{\xi}(m) \in \ker \omega_m$
if and only if $\ceV{\xi}(m) \in \ker \omega_m$.
\end{enumerate}
\end{cor}
\begin{pf}
To prove  (1), it suffices to show   that  if
 $\Vec{\xi} (m)+v \in \ker \omega_m$,
 where $\xi\in {A}_{m}$ and $v\in T_m P$,
 then both  $\Vec{\xi} (m)$ and $v$ belong to $\ker \omega_m$.
According to Proposition \ref{prop:2.2} (1),
for any $u\in  T_m P$, we have
$$\omega ( \Vec{\xi} (m), u)=\omega ( \Vec{\xi} (m)+v , u)=0. $$
On the other hand, for any $\eta \in {A}_{m}$,  we have
$\omega ( \Vec{\xi} (m), \ceV{\eta }(m))=0$ according to
Proposition \ref{prop:2.2} (3).
Thus it follows that $\Vec{\xi} (m) \in \ker \omega_m$, which
also implies that $v\in \ker \omega_m$.

(2)  Note that  $a (\xi )= \Vec{\xi}(m)-\ceV{\xi }(m)$.
Hence for any  $\eta \in {A}_{m}$, we have
$$\omega (a (\xi ), \Vec{\eta }(m) )=
\omega ( \Vec{\xi}(m)-\ceV{\xi }(m), \Vec{\eta }(m) )=
\omega ( \Vec{\xi}(m), \Vec{\eta }(m) )-\omega (\ceV{\xi }(m),  
\Vec{\eta }(m) )=0. $$
It thus follows that $a (\xi )\in \ker \omega_m $ since 
$\epsilon^* \omega =0$ according to Proposition \ref{prop:2.2} (1).

(3) follows from (2) since $a (\xi )= \Vec{\xi}(m)-\ceV{\xi }(m)$.
\end{pf}

\subsection{Quasi-symplectic groupoids}

Let us set

\begin{equation}
\label{eq:kerA}
 \ker \omega_m \cap  {A}_{m}=
\{\xi\in A_m | \Vec{\xi}(m)\in \ker \omega_m\}. 
\end{equation} 

Corollary \ref{cor:2.3} implies that the anchor
induces a well-defined map from $ \ker \omega_m \cap  {A}_{m}$
 to $\ker \omega_m \cap T_{m}P$.
Now we  are ready to introduce the non-degeneracy condition.

\begin{defn}
\label{def:non-deg}
A pre-quasi-symplectic  groupoid  $({\gm}\toto{P}, \omega +\Omega )$
is said to be quasi-symplectic if the following non-degeneracy 
 condition is satisfied: the anchor  
$$a: \ker \omega_m \cap  {A}_{m}\to 
\ker \omega_m \cap T_{m}P$$
 is an isomorphism.
\end{defn}  

Given 
a pre-quasi-symplectic   groupoid  $({\gm}\toto{P}, \omega +\Omega )$,
the two-form $\omega$ induces a well-defined  linear map:
$$\omega^b :T_m P\lon A^*_m, \ \ \ \la \omega^b (v), \xi \ra
=\omega (v, \Vec{\xi}(m) ), \ \ \forall v\in T_m P, \ \xi\in A_m .$$
 Indeed one  easily sees that
  $\omega^b$ induces a well-defined map:
\begin{eqnarray}
\label{eq:TA}
&&\phi: \frac{T_m P}{\ker \omega_m \cap T_{m}P}\lon
( \frac{A_m}{\ker \omega_m \cap 
A_m})^* , \nonumber \\
&& \la   \phi   [v], [\xi]\ra = \la \omega^b (v), \xi \ra 
=\omega (v, \Vec{\xi}(m)), 
\ \ \ \forall  v\in T_m P,  \ \ \xi \in A_m.
\end{eqnarray}

The following result plays an essential role in
understanding the non-degeneracy condition.

\begin{prop}
\label{pro:phi}
Assume that $({\gm}\toto{P}, \omega +\Omega )$
is a  pre-quasi-symplectic   groupoid. Then
 $\phi$  is a  linear isomorphism.
\end{prop}
\begin{pf}
   Assume that $\phi  [v]=0$ for $v\in  T_{m}P$.
Then $\omega (v, \Vec{\xi}(m))=0, \ \forall \xi \in A_m$,
which implies that $v \in \ker \omega_m $ since $\epsilon^* \omega =0$.
Hence $[v]=0$. So $\phi $ is injective.

Conversely, assume that $\xi \in A_m$ satisfies the property
that  $\la   \phi   [v], [\xi]\ra =0, \ \ \forall v\in T_m P$.
Hence $\omega (\Vec{\xi}(m), v)=0, \ \forall v\in T_m P$. This
implies that $\Vec{\xi}(m)\in \ker  \omega_m$.
Therefore $\xi \in  \ker \omega_m \cap  {A}_{m}$, or $[\xi]=0$.
This implies that $\phi  $ is surjective.
\end{pf}

An immediate consequence is the following result, which
 gives a useful way of characterizing
a quasi-symplectic groupoid.

\begin{prop}
\label{prop:eq}
A pre-quasi-symplectic   groupoid  $({\gm}\toto{P}, \omega +\Omega )$
is a    quasi-symplectic  groupoid if and only if
\begin{enumerate}
\item the anchor  $a: \ker \omega_m \cap  {A}_{m}\to
\ker \omega_m \cap T_{m}P$ is injective, and
\item   $\dim \gm =2 \dim P$.
\end{enumerate} 
\end{prop}
\begin{pf}
 By  Proposition \ref{pro:phi} and using dimension counting, 
we have
\begin{equation}
\label{eq:dim}
\dim (\ker \omega_m \cap  {A}_{m})-\dim (\ker \omega_m \cap T_{m}P)
=\dim \gm -2\dim P. 
\end{equation}

Assume that $({\gm}\toto{P}, \omega +\Omega )$ is a 
quasi-symplectic  groupoid. Eq. (\ref{eq:dim}) implies
that $\dim \gm =2 \dim P$. The converse is proved by working backwards.
\end{pf}

A special class of quasi-symplectic  groupoids are the 
so called {\em twisted symplectic  groupoids} \cite{CX},  which are
pre-quasi-symplectic   groupoids $({\gm}\toto{P}, \omega +\Omega )$
such that $\omega$ is honestly non-degenerate.   In particular, 
 symplectic groupoids \cite{Weinstein:87} are  always
quasi-symplectic.  In the next subsection,
we will discuss another  class of quasi-symplectic  groupoids motivated by the 
Lie group valued momentum map theory of
Alekseev--Malkin--Meinrenken \cite{AMM}.

\subsection{AMM quasi-symplectic groupoids}

First of all, let us fix some notations.
Assume that a Lie group $G$ acts smoothly  on a  manifold
 $M$  from the left.
 By a transformation groupoid, we mean the  groupoid
$G\times M\toto M$, where 
the source and target maps are given, respectively, by
 $s  (g, x)=gx, \  t  ({g}, {x})=x$,
$\forall (g, x)\in G\times M$,
 and the multiplication is
 $ ({g}_{1}, x )\cdot ({g}_{2} , y) =({g}_{1}{g}_{2},  y)$,
 where $ x={g}_2 y$. 

Let $G$ be a Lie group equipped with an ad-invariant
 non-degenerate symmetric bilinear form $(\cdot , \cdot )$.
 Consider the transformation groupoid
${G\times G}\toto {G}$, where $G$ acts on itself by  conjugation.
Following  \cite{AMM},  we denote by $\theta $ and $\bar{\theta}$
the left and right Maurer-Cartan forms on $G$ respectively, i.e.,
$\theta =g^{-1}dg$ and $\bar{\theta}=dg g^{-1}$.
Let $\Omega \in \Omega^{3}(G)$ denote the bi-invariant
3-form on $G$ corresponding to the Lie algebra 3-cocycle
$\frac{1}{12}(\cdot , [\cdot , \cdot ])\in \wedge^3 \frakg^*$:
\begin{equation}
\label{eq:chi}
\Omega =\frac{1}{12}(\theta , [\theta , \theta ])
=\frac{1}{12}(\bar{\theta}, [\bar{\theta} , \bar{\theta} ])
\end{equation}
and  $\omega \in \Omega^2 (G\times G )$  the two-form:

\begin{equation}
\label{eq:quasi}
\omega|_{(g, x)} =-\half [(Ad_{x} \pr_1^* \theta , \pr_1^* \theta )
+(\pr_1^* \theta , \pr_2^{*}(\theta +\bar{\theta} ))],
\end{equation}
where $(g, x)$ denotes the coordinate in $G\times G$, and
$\pr_1$ and $\pr_2: G\times G\to G$ are the natural projections.

\begin{prop}
Let $G$ be a Lie group equipped with an ad-invariant non-degenerate symmetric
bilinear form $(\cdot , \cdot )$. Then the transformation groupoid
$({G\times G}\toto {G}, \omega+\Omega )$ is a quasi-symplectic
groupoid, called the AMM quasi-symplectic groupoid.
\end{prop}
\begin{pf}
First, one needs to check that $\omega+\Omega $ is
a 3-cocycle. This can be done by a tedious computation,
and is left for the reader.

It remains to check the non-degeneracy condition, which is in fact
embedded in the proof  of  Proposition 3.2 \cite{AMM}.
For completeness, let us sketch a proof below.

The Lie algebroid $A$ of ${G\times G}\toto {G}$ is a transformation
Lie algebroid:
$\frakg \times G\to G$, where the anchor map $a: \frakg \times G\to TG$
is given by  $a (\xi, x )=r_x (\xi )-l_x (\xi )$, $ \forall
\xi \in \frakg$.
Therefore  $a (\xi, x )=0$ if and only if $Ad_x \xi =\xi$.
On the other hand, for any $\xi \in \frakg$ being identified
with an element in $A_x$, we have 
$\Vec{\xi}|_{(1, x)}=(\xi, 0)\in T^{t}_{(1, x)} (G\times G)$.
For any $\delta_x \in T_x G$, let $\delta_{(1, x)}=(0, \delta_x )
\in T_{(1, x)} (G\times G)$. Clearly $\delta_{(1, x)}$ is
a tangent vector to the unit space.

It  follows from Eq. (\ref{eq:quasi}) that
$$\omega (\Vec{\xi}|_{(1, x)},  \delta_{(1, x)} )=
\omega ((\xi, 0), (0, \delta_x ))=
-\half  \delta_x \per (\xi,  \theta +\bar{\theta } ).$$
Therefore we have
$\epsilon^* (\Vec{\xi}|_{(1, x)} \per \omega )=\half  (\xi, 
 \theta +\bar{\theta } )$. Hence, $\Vec{\xi}|_{(1, x)} \per \omega=0$
 if and only if $(Ad_x +1)\xi =0$. This implies that
   $a: \ker \omega_x \cap  {A}_{x}\to
\ker \omega_x \cap T_{x}G$ is injective. Therefore it follows 
from Proposition  \ref{prop:eq} that $({G\times G}\toto {G}, \omega+\Omega )$
 is indeed a quasi-symplectic
groupoid.
\end{pf}


\begin{numrmk}
From the above proposition, we see that $[\omega+\Omega ]$
defines a class in the equivariant cohomology
$H^{3}_{G}(G)$.
 When $G$ is a   compact simple Lie group with the  basic
form $(\cdot , \cdot )$,
$[\omega+\Omega ]$
is a generator of   $H^{3}_{G}(G)$.
 In Cartan model, it corresponds
to the class defined by the $d_{G}$-closed equivariant 3-form
$\chi_{G}(\xi )=\Omega-\half ( \theta +\bar{\theta}, \xi ):
\ \frakg \lon \Omega^* (G), \ \forall \xi\in \frakg$ (see \cite{BCWZ, M}).
\end{numrmk}

\section{Hamiltonian $\gm$-spaces}

\subsection{Definitions and properties}

In this subsection, we introduce the notion of Hamiltonian $\gm$-spaces
for a quasi-symplectic groupoid $\gm\toto P$, which generalizes 
the usual notion of Hamiltonian spaces of symplectic groupoids
in the sense of Mikami-Weinstein \cite{MW}.

First, we need the following:

\begin{defn}
\label{def:gm-space}
Given a quasi-symplectic  groupoid  $({\gm}\toto{P}, \omega +\Omega )$, 
let $J: X\to P$ be a left $\gm$-space, i.e., $\gm$ acts on  $X$ from
the left. By a  compatible  two-form on $X$, we mean a 
 two-form $\omega_X \in \Omega^2 (X)$  satisfying
\begin{enumerate}
\item $d \omega_X =J^* \Omega $; and
\item the graph of the action $\Lambda=\{(r, x, rx)|t(r) =J(x)\}
 \subset \gm \times X \times X$
is isotropic with respect to the two-form
 $(\omega,  \omega_X, - \overline{\omega_X  } )$.
\end{enumerate}
Then $(X\stackrel{J}{\to}P, \omega_X)$ is called a  pre-Hamiltonian
$\gm$-space. 
\end{defn} 

In the sequel,  we  simply refer to the second condition as  to 
``the graph of the action $\Lambda \subset \gm \times X \times \ol{X}$
is isotropic", where the bar on the last factor $X$ indicates
that the  opposite two-form is used.

To illustrate the intrinsic meaning of the above
  compatibility condition, let us elaborate it in terms of groupoids.
Let $Q :=\gm \times_P X \toto X$ denote the transformation groupoid
corresponding to the $\gm$-action,
 and, by abuse of notation,  $J :Q\to \gm$ the natural projection.
It is simple to see that
%

\begin{equation}
\xymatrix{
Q \ar@<-.5ex>[d]\ar@<.5ex>[d]\ar[r]^{J}
&
\gm \ar@<-.5ex>[d]\ar@<.5ex>[d]\\
X\ar[r]^J & P}
\end{equation}
is a Lie groupoid homomorphism. Therefore it induces a map, i.e.,
 the pull-back map, on the level of de-Rham complex
$$J^* :\ \ \  \Omega^{\com}(\gm\lcom )\to \Omega^{\com}(Q\lcom ). $$

\begin{prop}
\label{pro:3.2}
Let $({\gm}\toto{P}, \omega +\Omega )$ be a quasi-symplectic  groupoid
and $J: X\to P$ a left $\gm$-space. Then
 $\omega_X\in \Omega^2 (X)$ is
a compatible two-form if and only if
\begin{equation}
\label{eq:pr}
J^* ( \omega +\Omega )=\delta \omega_X.
\end{equation} 
\end{prop}
\begin{pf}
Note that
$$\delta \omega_X=(s^*\omega_X-t^*\omega_X)+ d\omega_X, $$
where $s, t:  \gm \times_P X \to X$ are the source and
target maps of the groupoid $\gm \times_P X\toto X$. 
So  Eq. (\ref{eq:pr}) is equivalent to 
$$ s^*\omega_X-t^*\omega_X =J^*  \omega,  \ \ \ \mbox{and }
 d\omega_X=J^* \Omega.$$
It is simple to see that the first  equation above
 is equivalent to that 
the graph of the action $\Lambda \subset \gm \times X \times \ol{X}$
is isotropic by using  the source and target maps 
 $s(r, x)=r\cdot x$ and $t(r, x)= x$, $\forall (r, x)\in \gm \times_P X$.
\end{pf} 

\begin{numrmk}
As a consequence, $J^* : H^3_{DR} (\gm\lcom)\to H^3_{DR} (Q\lcom)$
maps $[\omega +\Omega]$ into zero. When  $[\omega +\Omega]$ is of integral
class, it defines an $S^1$-gerbe over the stack $\XX_\gm$ corresponding
to the groupoid $\gm \toto P$, the above
proposition implies that the pull-back $S^1$-gerbe  on $\XX_Q$
is always trivial. 

If $\gm$ is the symplectic groupoid $T^* G\toto \frakg^*$,
 $Q$ can be identified with the
transformation groupoid $G \times X\toto X$ and the 
groupoid homomorphism $J: Q\ (\cong G \times X)\to 
\gm \ (\cong G \times \frakg^* )$ is simply
 $id \times J$.  In this case, $H^3_{DR} (\gm\lcom) \cong H^3_G (\frakg^* )$
and $ H^3_{DR} (Q\lcom) \cong H^3_G (X)$. 
In Cartan model, Eq. (\ref{eq:pr}) is equivalent to
\begin{equation*}
d_G  \omega_X= 
 J^* \chi_G(\xi).
\end{equation*}
Here $\chi_G\in \Omega^3_G(\frakg^*)$ is the
 equivariant  closed 3-form defined
as  $\chi_G(\xi)=-d\la a,\xi\ra$, where $a:\frakg^*\to\frakg^*$
 is the identity map. Similarly, if  $\gm$ is the AMM quasi-symplectic
groupoid $G\times G\toto G$, $Q$ is isomorphic
to the transformation groupoid $G\times X\toto X$.
Then
the relevant $d_G$-closed equivariant 3-form
$\chi_G\in\Omega^3_G(G)$ is
\begin{equation*}
\chi_G(\xi)=\Omega -\frac{1}{2}(\theta+\bar{\theta},\xi ) .
\end{equation*}
See Remark 2.1 of \cite{AMM}.

Note that in the first case,  $\chi_G \in \Omega^3_G(\frakg^* )$
 defines a zero class in $H^3_G (\frakg^* )$, while in
the case of the AMM quasi-symplectic groupoid,
$\chi_G \in\Omega^3_G(G)$ defines a non-zero class
in $H^3_G (G)$. This fact is the key ingredient  for
explaining the difference of their quantization theories, while
in the latter  case, gerbes are inevitable in the construction \cite{LX}.  
\end{numrmk}

As is well known, a Lie groupoid action induces a Lie algebroid
action, called the infinitesimal action, which can be described
as follows. 
For any $x\in X$ and any $\xi \in A_m$,  where $J(x)=m$,
let $\gamma (t)$ be a  path in the $t$-fiber $t^{-1}(m)$ of 
$\gm$ through
the point $m$ such that $\dot{\gamma }(0)=\Vec{\xi} (m) $,
and define $\hat{\xi}(x) \in T_x X$  to be the tangent vector  corresponding
to the curve $\gamma (t)\cdot x$ through the point $x$.
In this way  one obtains a linear map 
$$A_m\lon T_x X, \ \ \ \xi \to \hat{\xi}(x)$$
called the infinitesimal action. In particular, this action induces
a Lie algebra homomorphism $\gm (A)\to {\mathfrak{X}} (X)$.
 One also easily checks that
$$a(\xi )=J_{*}\hat{\xi}(x), \ \ \ \forall \xi \in A_m .$$

The following lemma follows easily from the compatibility condition
in Definition \ref{def:gm-space} (2).

\begin{lem}
\label{lem:mom}
Let $(\gm \toto P, \omega+\Omega )$ be a quasi-symplectic groupoid.
If a $\gm$-space $J: X\to P$ equipped with a two-form
$\omega_X$ satisfies the compatibility condition
in Definition \ref{def:gm-space} (2),
then for any $x\in X$ such that $J(x)=m$ and any $\xi \in A_m$, 
we have
\begin{equation}
\label{eq:mom}
J^*  \epsilon^* (\Vec{\xi} (m) \per \omega )=
\hat{\xi}(x)\per \omega_X .
\end{equation} 
\end{lem}  
\begin{pf}
It is simple to see that
for any $\xi \in A_m$,    $(\Vec{\xi} (m) , 0 , \hat{\xi}(x))$
is tangent to $\Lambda$. On the other hand,
$\forall \delta_x \in T_x X$, $(J_* \delta_x , \delta_x ,\delta_x )$
is also tangent to $\Lambda$. Thus it follows that
$$\omega (\Vec{\xi} (m) , J_* \delta_x )-\omega_X (\hat{\xi}(x), \delta_x )
=0.$$
Eq. (\ref{eq:mom}) thus follows immediately. 
\end{pf}                 
 
From this lemma, one easily sees that if $\Vec{\xi} (m)  \in
\ker \omega$,  then $\hat{\xi}(x)$  automatically
 belongs to the kernel of $\omega_X $.
As in  \cite{AMM}, we impose the following minimal non-degeneracy
 condition.

\begin{defn}
Let $({\gm}\toto{P}, \omega +\Omega )$ be a quasi-symplectic  groupoid.
A Hamiltonian $\gm$-space is  a  left $\gm$-space $X\to P$ 
equipped with a compatible two-form $\omega_X$ such that $\forall x\in X$,
\begin{equation}
\ker \omega_X |_x =\{ \hat{\xi}(x)| \xi \in A_{J(x)} \ \mbox{ such that }
\Vec{\xi} (J(x))  \in \ker \omega\}.
\end{equation}
\end{defn}

For any $x\in X$, by 
   $A^{x}_{x}$, we denote the linear subspace of $A_{J(x)}$
consisting of those vectors $ \xi \in A_{J(x)}$ such that
$\hat{\xi}(x)=0$.

\begin{lem}
\label{lem:Ga-space}
Assume that $({\gm}\toto{P}, \omega +\Omega )$ is  a quasi-symplectic
groupoid and  $J: X\to P$ is  a  $\gm$-space equipped
with  a compatible two-form $ \omega_X $.
Then
\begin{enumerate}
\item  $\dim J_{*}(T_x X)\leq \rank A-\dim A^{x}_{x}$;
\item if moreover $( X\stackrel{J}{\to} P,  \omega_X )$ 
is an  Hamiltonian $\gm$-space, then

{\bf a.} $\ker \omega_{J(x)}\cap A_{J(x)} \to \ker \omega_X |_x$,
 $ \xi \to \hat{\xi}(x)$ is an isomorphism; and

{\bf b.}   $ \ker J_{*} \cap \ker \omega_X^b=0$.
\end{enumerate}
\end{lem}
\begin{pf}
   (1) $\forall \delta_x\in T_x X$ and $\xi \in A^{x}_{x}$, we have
$$\la \phi [J_* \delta_x ],  [\xi ]\ra
= \omega ( J_* \delta_x, \Vec{\xi}(J(x) ))
= -\omega_X (\hat{\xi }(x), \delta_x)=0,$$
where $\phi$   is the linear isomorphism defined by Eq. (\ref{eq:TA}).
This implies  that  
$$\phi [ \pr_1 (J_* (T_x X)) ]\subseteq  (\pr_2 A^{x}_{x} )^{\perp},$$
where 
$$\pr_1 : T_{J(x)} P\to \frac{T_{J(x)}  P}{\ker \omega_{J(x)} \cap T_{J(x)}P}$$
 and
$$\pr_2:  A_{J(x)}  \to \frac{A_{J(x)} }{\ker \omega_{J(x)}  \cap  A_{J(x)} }$$
 are   projections.

Secondly, we note that $\pr_2$ is injective when being restricted
to $ A^{x}_{x}$. To see this, we only  need to show that
$A^{x}_{x} \cap  (\ker \omega_{J(x)} \cap  A_{J(x)}) =0$.
Assume that $\xi \in A^{x}_{x} \cap  (\ker \omega_{J(x)} \cap  A_{J(x)})$.
Then we have $\hat{\xi }(x)=0$ and $ \Vec{\xi}(\mmm )\per \omega =0$.
Hence $ a(\xi )=J_* \hat{\xi }(x)=0$, which implies that
$\xi =0$ by Definition \ref{def:non-deg}. 
 As a consequence, we have $\dim (\pr_2 A^{x}_{x})=\dim
A^{x}_{x}$. Hence,

\be
&&\dim J_{*}(T_x X) -\dim (\ker \omega_{\mmm}\cap T_{\mmm} P)\\
&\leq &\dim \pr_1 (J_{*}(T_x X)) \ \ \ \mbox{(since $\phi$ is a linear
isomorphism)}\\
&=&\dim \phi [\pr_1 (J_{*}(T_x X)) ]\\
&\leq &\dim  (\pr_2 A^x_{x} )^{\perp}\\
&= & [\rank A -\dim   (\ker \omega_\mmm \cap A_\mmm )]-\dim (\pr_2 A^x_{x} )\\
&= & [\rank A -\dim   (\ker \omega_\mmm \cap A_\mmm )]-\dim  A^x_{x}.
\ee
Thus (1) follows immediately since $\gm$ is a quasi-symplectic groupoid.

(2) (a). By the minimal  non-degeneracy assumption,
we  know that the map 
$$  \ker \omega_{\mmm}\cap A_{\mmm} \to \ker \omega_X|_x,\ \   \xi
 \to \hat{\xi}(x),$$
is surjective. To show that it is injective,
 assume that $\xi \in \ker \omega_{\mmm}\cap A_{\mmm}$
 such that $\hat{\xi}(x) =0$.
Then $a(\xi )=J_{*}\hat{\xi}(x) =0$. Since $\omega $ is
 non-degenerate in the sense of
Definition \ref{def:non-deg}, we have $\xi =0$.

 (b).  Assume that $\delta_x \in \ker J_* \cap  \ker \omega_X^b $.
Since $J: X\lon P$ is an Hamiltonian $\gm$-space,
by assumption, we have $\delta_x =\hat{\xi}(x)$, where   $\xi\in A_{J(x)}$
such that $\Vec{\xi} (\mmm)   \per \omega =0$.
Hence $a (\xi )=J_* \hat{\xi}(x) =J_* \delta_x =0$,
and therefore $\xi=0$ since $\gm$ is a quasi-symplectic groupoid.

This completes the proof. 
\end{pf}

For a subspace $V\subseteq T_x X$, by $V^{\omega_X}$ we denote
its $\omega_X$-orthogonal subspace of $V$.
As a consequence, we have the following proposition which 
plays a key role in our reduction theory.

\begin{prop}
Assume that  $({\gm}\toto{P}, \omega +\Omega )$
is a  quasi-symplectic  groupoid, and
  $(X\stackrel{J}{\to}P,  \omega_X )$
 an  Hamiltonian $\gm$-space. Then  
\begin{equation}
\label{eq:key}
(\ker J_* )^{\omega_X}=\{\hat{\xi}(x)|\forall \ \xi\in A_{J(x)}\}.
\end{equation} 
\end{prop}
\begin{pf}
  It is simple to see that $(\ker J_* )^{\omega_X}=
[\omega_X^b (\ker J_* )]^{\perp}$. Therefore it follows that
\be
&&\dim (\ker J_* )^{\omega_X}\\
&=&\dim X-\dim [\omega_X^b (\ker J_* )]\ \ \ \ \mbox{(since $\omega_X^b$
is injective when being restricted to $\ker J_* $)}\\
&=&\dim X-\dim\ker J_* \\
&=&\dim J_* (T_x X  )\ \ \ \ \mbox{(by Lemma \ref{lem:Ga-space})}\\
&\leq &\rank A-\dim A_{x}^x\\
&=&\dim \{\hat{\xi}(x)|\forall \ \xi\in A_{J(x)}\}.
\ee
On the other hand, clearly we have 
$$\{\hat{\xi}(x)|\forall \ \xi\in A_{J(x)}\}\subseteq (\ker J_* )^{\omega_X}$$
according to Eq. (\ref{eq:mom}). 
Thus Eq. (\ref{eq:key})  follows immediately. 
\end{pf}


\subsection{Two fundamental examples}

Below we study two fundamental 
examples of Hamiltonian $\gm$-spaces, which  are naturally   associated 
to a  quasi-symplectic groupoid.

\begin{prop}
\label{prop:example}
Assume that  $({\gm}\toto{P}, \omega +\Omega )$ is   a 
quasi-symplectic groupoid. Then
\begin{enumerate}
\item $J: \gm \to P\times {P}$ is  an Hamiltonian 
$\gm \times \overline{\gm }$-space, where $J(r)=(s(r), t(r)), \ \forall
r\in \gm$,  and the action is defined by
$$ (r_1 , r_2 )\cdot x=r_1  x r_2^{-1},   \ \ \ t(r_1 )=s (x), \ t(x)=t (r_2 ).$$
\item Given any orbit $\O\subset P$, there is a natural two-form
 $\omega_{\O}\in \Omega^2 (\O)$ so that the natural
inclusion  $i: \O\to P$
 defines an Hamiltonian $\gm $-space under the natural $\gm $-action.
\end{enumerate}  
\end{prop}
\begin{pf} 
(1) It is clear,  from definition,  that $d \omega =J^*\Omega$.
To check  the second  compatibility condition of Definition 
\ref{def:gm-space}, it suffices to show that
$$\{(r_1 , r_2 , x, r_1  x r_2^{-1} )| \ t(r_1 )=s (x), \ t(x)=t (r_2 )\}
\subset \gm \times \overline{\gm }\times \gm \times \overline{\gm }$$
is isotropic. This can be proved using the multiplicativity
assumption on $\omega$, i.e., $\partial \omega =0$,
 as in \cite{Weinstein:affinoid}.
To check the minimal non-degeneracy condition, note  that for any
$\xi, \eta \in \gm (A)$, the vector field on $\gm$ generated by
the infinitesimal action of
$(\xi, \eta)$ is given by $\Vec{\xi}(x)-\ceV{\eta}(x)$.
Next, note that for any $\delta_x\in T_x \gm, \ \xi\in \gm (A)$,
we have

\begin{equation}
\label{eq:left-right}
\omega (\ceV{\xi}(x), \delta_x )
=\omega (\ceV{\xi} (t(x)), t_{*}\delta_x ), 
\ \ \ \omega (\Vec{\xi}(x), \delta_x )
=\omega (\Vec{\xi} (s(x)), s_{*}\delta_x ).
\end{equation} 

These  equations follow  essentially  from Eq. (\ref{eq:mom})
  since  $s: \gm \to P$ 
 equipped with the natural left $\gm$-action
(or $t:  \gm \to P$  with the left $\gm$-action: $r\cdot x=x r^{-1}$,
 respectively)
satisfies the hypothesis of Lemma  \ref{lem:mom}.

Now assume that $\delta_x\in \ker \omega$.
 Then $t_{*}\delta_x \in \ker \omega$ by Eq. (\ref{eq:left-right}), 
since $P$ is isotropic. By the non-degeneracy assumption,
we have $t_{*}\delta_x =a (\eta )$ for some $\eta
\in A|_{t(x)}$ such that  $\Vec{\eta } (t(x)) \in
\ker \omega$.  Hence  $\ceV{\eta}(t(x))\in
\ker \omega$ by Corollary   \ref{cor:2.3} (3), which in
turn implies that $\ceV{\eta}(x)\in \ker \omega$
 according to Eq. (\ref{eq:left-right}).
  Let $\delta'_x =\delta_x +\ceV{\eta}(x)$. Then,
$$t_{*} \delta'_x =t_{*} \delta_x +t_{*}\ceV{\eta}(x)=t_{*} \delta_x -a (\eta )
=0.$$
Also we know that $\delta'_x \in \ker \omega $. Therefore
one can write $\delta'_x =\Vec{\xi}(x)$ where $\xi \in
A_{s(x)}$ such that  $\Vec{\xi}( s(x) )\in 
\ker \omega $.
We thus have  proved that  $\delta_x =\Vec{\xi}(x)-\ceV{\eta}(x)$,
where $\Vec{\eta} (t (x) )\in \ker \omega$ and $ \Vec{\xi}( s(x))\in 
\ker \omega $.

(2) Let $\O\subset P$ be the groupoid
 orbit through  the point $m_0\in P$.
It is standard that $t^{-1}(m_0)\stackrel{s}{\to}\O$ is a 
$\gm_{m_0}^{m_0}$-principal bundle, where $\gm_{m_0}^{m_0}$
 denotes the isotropy group at $m_0$.
 From the multiplicativity assumption on  $\omega$,
it is simple to see that  $\omega|_{t^{-1}(m_0)}$, the pull-back
of $\omega$ to the $t$-fiber $t^{-1}(m_0)$, is indeed basic
with respect to the $\gm_{m_0}^{m_0}$-action.
 Hence it descends to a two-form
 $\omega_\O$ on $\O$. That is, $\omega|_{t^{-1}(m_0)}= s^* \omega_\O$.
It thus follows that 
$$s^* d\omega_\O =(s^*\Omega -t^* \Omega ) |_{t^{-1}(m_0)},$$
which implies that $d\omega_\O =i^* \Omega $.
It is also clear that  the two-form $\omega_\O$ is compatible with
the groupoid $\gm$-action since $\omega$ is multiplicative.
To show the minimal  non-degeneracy condition, assume that  
$x\in t^{-1}(m_0)$ is an arbitrary point, and $\delta_x \in T_x t^{-1}(m_0)$ 
such that $[\delta_x]=s_{*} \delta_x \in  \ker \omega_{\O}|_m$, where
$m=s(x)$.
By definition,
$\omega  ( \delta_x, \delta'_x )=0, \ \forall \delta'_x\in
 T_x t^{-1}(m_0)$.
It thus follows that
 $\omega  (r_{x^{-1}}\delta_x, r_{x^{-1}} \delta'_x )=0$.
Let $\xi, \eta \in A_{m}$
 such that $r_{x^{-1}}\delta_x =\Vec{\xi} (m)$
and $r_{x^{-1}}\delta'_x =\Vec{\eta} (m)$. Thus
we have $\omega (\Vec{\xi} (m), \Vec{\eta} (m) )=0, \ \forall \eta\in A_{m}$.
Therefore 
$$\omega (a(\xi ), \Vec{\eta} (m) )=\omega  (\Vec{\xi} (m)-\ceV{\xi} (m), \ \Vec{\eta} (m) )=\omega(\Vec{\xi} (m) , \Vec{\eta} (m) )=0, \ \forall \eta
\in A_m. $$
It thus follows that 
$a(\xi )\in \ker \omega$ since $\omega (a(\xi ), T_{m}P)=0$.
 That is, $a(\xi )\in \ker \omega_m \cap T_m P$.
 By the non-degeneracy assumption on
$\omega $ (see Definition \ref{def:non-deg}),
 we deduce  that there  exists
$\xi_1\in A_m$ such that
$\Vec{\xi_1} (m)\in \ker \omega$ and $a(\xi_1 )=a(\xi )$. So
$\xi -\xi_1$ belongs to the isotropy Lie algebra at $m$. 
As a result,  it follows that
the minimal non-degeneracy  condition is  indeed satisfied since
$ [\delta_x ]=s_{*} \delta_x =\hat{\xi }(m) =\hat{\xi_1}(m)$.
\end{pf}

\subsection{Examples of  Hamiltonian $\gm$-spaces}

In this subsection, we list various examples of momentum
maps appeared in the literature, which
 can be considered as special cases of our
Hamiltonian $\gm$-spaces.
In fact, our definition is a natural generalization
of Hamiltonian $\gm$-spaces of a  symplectic
groupoid of Mikami--Weinstein \cite{MW}, which include the
usual Hamiltonian  momentum maps and Lu's momentum maps of Poisson group
actions as special cases.

\begin{numex}
Consider the symplectic groupoid $(T^*G \toto \Gg^*, \omega )$,
where $\omega$ is the standard cotangent  symplectic structure.
Then its Hamiltonian spaces are exactly the Hamiltonian
$G$-spaces $J: X\to \Gg^*$ in the ordinary sense.
\end{numex}       

\begin{numex}

When $P=G^{*}$, the dual of a simply connected
  complete Poisson Lie group $G$, its
  symplectic groupoid $\Gamma$ is a transformation
groupoid: $G\times G^* \toto G^*$, where $G$ acts on $G^*$
by left dressing action \cite{LuW:1989}.  In this case, 
Hamiltonian $\Gamma$-spaces can be
  described in terms of the so-called Poisson $G$-spaces.
  Recall that a  symplectic (or more generally a Poisson)
  manifold $X$ with a left $G$-action is called a Poisson $G$-space
  if the action map $\overline{G}\times X\rightarrow X$ is a
 Poisson map.  A Poisson morphism $J:X\rightarrow G^*$ is said to 
be a momentum
 map for the Poisson $G$-space \cite{Lu:1991}, if
  \begin{equation}
 \label{eq:momentum-group}
 X\in \frakg  \mapsto  -\pi^{\#}_{X}(J^{*}(X^{r}))\in \XX (X)
 \end{equation}
 is the
  infinitesimal generator of the $G$-action, where $X^{r}$ denotes
  the right-invariant one-form on $G^{*}$ with value $X\in\frakg^*$ at
  the identity, and $\pi_X$ is the Poisson tensor on $X$.
An explicit relation between Hamiltonian $\gm$-spaces 
and Poisson $G$-spaces can be established as 
follows \cite{WeinsteinX:yang}.
 If $J:X\to  G^{*}$ is an Hamiltonian 
   $\Gamma$-space, then $X$ is a Poisson $G$-space with the
   action:
   \begin{equation} \label{eq:action}
   gx=(g,J(x))\cdot x, \end{equation}
   for any $g\in G$ and $x\in X$, where
   $(g,J(x))$ is considered as an element in $\Gamma=
   G\times G^{*}$ and the dot on the right hand side refers to the
groupoid
 $\Gamma$-action on $X$. Then $J$ is the momentum
   map of the induced Poisson $G$-action, in the sense 
of Lu \cite{Lu:1991}.
     Conversely, if a symplectic manifold $X$ is a Poisson
   $G$-space with a momentum mapping $J:X\to G^{*}$, Eq.
   (\ref{eq:action}) defines an Hamiltonian  $\gm$-space.
\end{numex}

\begin{numex}
Let $(\cdot, \cdot )$ be an ad-invariant 
 non-degenerate symmetric bilinear
form on $\frakg$. It is well-known that
$(\cdot, \cdot )$ induces a Lie algebra 2-cocycle 
$\lambda \in \wedge^2 (L\frakg^* )$ on the
loop Lie algebra defined by \cite{PS}:
\begin{equation}
\label{eq:lambda}
\lambda (X , Y)=\frac{1}{2\pi}\int_{0}^{2\pi}(X(s), Y'(s))ds, \  \forall
X(s) , Y(s)\in L\frakg,
\end{equation}
and therefore  defines an affine Poisson structure on $L\frakg $.
Its symplectic groupoid $\gm$ can be identified with
 the transformation groupoid
$LG \times L\frakg \toto  L\frakg $, where $LG$ acts on $L\frakg$
by the gauge transformation \cite{BXZ}: 
\begin{equation}
\label{eq:gauge}
g\cdot \xi = Ad_{g}\xi +  g'  g^{-1}, 
\ \ \forall g\in LG, \ \xi\in L\frakg.
\end{equation}  
This is  the standard  gauge transformation when
$L\frakg$  is  identified with the space of connections on  the
trivial $G$-bundle over the unit circle $S^1$.
    The symplectic structure on $LG \times L\frakg$
can be obtained as follows.                
By $\widetilde{L\frakg}$ we denote the corresponding  Lie algebra
central extension.
Assume that $\lambda $ satisfies  the integrality  condition (i.e.,
the corresponding closed two-form $\omega_{LG} \in  \Omega^2 (LG)^{LG}$
is of integer class). It defines  a loop group
central extension $S^1 \lon \widetilde{LG} \stackrel{\pi}{\lon} LG$.
Consider  $\tilde{\pi}:
\widetilde{LG} \times L\frakg \to {LG} \times L\frakg$, where
$\tilde{\pi}=\pi\times id$. Let $i$ denote the embedding
$\widetilde{LG} \times L\frakg\cong \widetilde{LG} \times (L\frakg \times
\{1\})\subset  \widetilde{LG} \times \widetilde{L\frakg} \cong T^* \widetilde{LG}$. Then
$$\tilde{\pi}^* \omega_{LG \times L\frakg}
=i^* \omega_{ T^* \widetilde{LG}}.$$

In this case, the corresponding Hamiltonian $\gm$-spaces are exactly
Hamiltonian loop group  spaces  studied extensively
by Meinrenken--Woodward \cite{MeinW, MeinW1, MeinW2}.
\end{numex}

\begin{numex}
Let $\gm$ be the AMM
quasi-symplectic groupoid  $(G\times G\toto G, \omega +\Omega )$.
It is simple to see that Hamiltonian $\gm$-spaces
 correspond exactly to
 quasi-Hamiltonian $G$ spaces with a group valued momentum map $J:X\to G$
 in the sense of \cite{AMM},  namely those $G$-spaces $X$ equipped
with a  $G$-invariant two-form 
 $\omega_X\in\Omega(X)^G$ and an equivariant map $J\in C^\infty(X,G)^G$
 such that: \begin{enumerate} 
\item[(B1)]\label{B1} The differential  of $\omega_X$ is given by: 
\begin{equation*}  \label{Gclosed} 
d\omega_X=  J^*\Omega.
 \end{equation*}                     
\item[(B2)]\label{B2}
The map $J$ satisfies
\begin{equation*}
\hat{\xi}\per \omega_X =\half J^*  (\xi, \theta +\bar{\theta }),
\forall \xi \in\frakg.
\end{equation*}
\item[(B3)]\label{B3}
At each $x\in X$, the kernel of $\omega_X$ is given by
\begin{equation*} 
\ker \omega_X|_x= \{ \hat{\xi}(x)|\ \xi\in\ker(\Ad_{J(x)}+1) \} .
\end{equation*}                      \end{enumerate}
\end{numex}

\subsection{Hamiltonian bimodules}

A useful way to study Hamiltonian $\gm$-spaces is via the Hamiltonian
bimodules.

\begin{defn}
Given quasi-symplectic groupoids $({G}\toto G_0 , \omega_G +\Omega_G )$
and $({H}\toto H_0 , \omega_H +\Omega_H )$, an Hamiltonian
$G$-$H$-bimodule is a 
manifold $X$  equipped with a two-form $\omega_X \in \Omega^2 (X)$ 
such that
\begin{enumerate} 
\item $G_{0}\stackrel{\rho }{\leftarrow} X
\stackrel{\sigma}{\rightarrow} H_{0}$ is  a left $G$-space  and
a right $H$-space, and   the two actions commute;
\item $X\stackrel{\rho \times \sigma}{\lon} G_{0}\times H_{0}$
is an Hamiltonian $G\times \overline{H}$-space, where the action is
given by
$(g, h)\cdot x=gxh^{-1}, \ \ \forall g\in G, \ h\in H, \ x\in X$ such that
$t(g)=\rho (x)$ and $t(h) =\sigma (x)$.
\end{enumerate}
\end{defn}

In particular, an Hamiltonian $\gm$-space can be considered
as  an  Hamiltonian $\gm$-$\cdot$-bimodule,
where $\cdot$ denotes the trivial quasi-symplectic  groupoid $\cdot\toto \cdot$.

Given an  Hamiltonian $G$-$H$-bimodule $G_{0}\stackrel{\rho }{\leftarrow} X
\stackrel{\sigma}{\rightarrow} H_{0}$,  let $Q\toto X$ be the
transformation groupoid  
$$Q:=(G\times H)\times_{(G_0\times H_0 )} X\toto X. $$
Then the  natural projections $\pr_1: Q\to G$ and  
$\pr_2: Q\to H$ are groupoid homomorphisms.
As an immediate consequence of Proposition \ref{pro:3.2},
 we have the following

\begin{prop}
\label{pro:pull}
If  $({G}\toto G_0 , \omega_G +\Omega_G )$
and $({H}\toto H_0 , \omega_H +\Omega_H )$
are quasi-symplectic groupoids, and 
 $G_{0}\stackrel{\rho }{\leftarrow} X
\stackrel{\sigma}{\rightarrow} H_{0}$ is an Hamiltonian 
$G$-$H$-bimodule, then
$$\pr_1^* (  \omega_G +\Omega_G )- \pr_2^*  (  \omega_H +\Omega_H )
=\delta \omega_X . $$
\end{prop}
Therefore, on the level of cohomology, we have 
$$\pr_1^* [\omega_G +\Omega_G]= \pr_2^*  [  \omega_H +\Omega_H ],$$
where $\pr_1^* : H^3_{DR}(G\lcom )\to H^3_{DR}(Q\lcom )$
and $\pr_2^* : H^3_{DR}(H\lcom )\to H^3_{DR}(Q\lcom )$
are the homomorphisms of cohomology groups induced
by the groupoid homomorphisms 
 $\pr_1: Q\to G$ and $\pr_2: Q\to H$, respectively.

Let $(G\toto G_0, \omega_G +\Omega_G)$, $(H\toto H_0, \omega_H +\Omega_H)$,
  and $(K\toto K_0, \omega_K +\Omega_K )$
be quasi-symplectic groupoids.
Assume that  $G_{0}\stackrel{\rho_1 }{\leftarrow} X
\stackrel{\sigma_1}{\rightarrow} H_{0}$ is an
 Hamiltonian
 $G$-$H$-\bimodule, and 
$H_{0}\stackrel{\rho_2 }{\leftarrow} Y
\stackrel{\sigma_2 }{\rightarrow} K_{0}$ 
 an  Hamiltonian $H$-$K$-bimodule. Moreover, we assume that
the fiber product  $X\times_{H_0}Y$ is
a manifold (for instance, this is true if $\sigma_1   \times \rho_2 :
X\times Y\to H_{0}\times  H_{0}$ is transversal to the
diagonal) and the diagonal $H$-action on $X\times_{H_0}Y$,
$h \cdot (x, y) =(x\cdot h^{-1}, h \cdot y)$, 
is free and proper so that the quotient space is a 
smooth manifold, which is denoted by $X\times_{H}Y$. That is
$$
X\times_{H}Y := \frac{X\times_{H_{0}}Y}{H}.
$$

Let $\rho_{3}: X\times_{H}Y\to G_{0}$ and $\sigma_{3}: 
X\times_{H}Y \to K_{0}$ be the maps given by
$\rho_{3}([x,y]) = \rho_{1}(x)$ and $\sigma_{3}([x,y])=\sigma_{2}(y)$,
respectively.
Define a left $G$-action and a right $K$-action on $X\times_{H}Y$ by
\begin{equation}
\label{eq:GK}
g\cdot [x,y]=[g\cdot x, y] \quad \text{and} \quad [x,y]\cdot k=[x, y\cdot k],
\end{equation}
whenever they are defined. 
It is clear that $ G_{0}\stackrel{\rho_3 }{\leftarrow}X\times_{H}Y
\stackrel{\sigma_3 }{\rightarrow} K_{0}$ becomes a 
left $G$- and right $K$-space, and that these two actions
commute with each other. 

To continue our discussion, we need to make  a  technical
assumption.

\begin{defn}
  We say that  two smooth maps  $\tau_i: X_i\to M$, $i=1, 2$,
are clean,  if
\begin{enumerate}
\item the fiber product $X_1 \times_{M} X_2$ is a smooth manifold; and
\item  for any $(x_1, x_2 )\in X_1 \times_{M} X_2$, $f_* T_{(x_1, x_2)}
(X_1 \times_{M} X_2)$ is equal to either $\tau_{1*}T_{x_1} X_1$ or 
$\tau_{2*}T_{x_2} X_2$, where $f: X_1 \times_{M} X_2\to M$ is defined 
as $f(x_1, x_2 )=\tau_1 (x_1)=\tau_2 (x_2 )$.
\end{enumerate}
\end{defn} 

For instance, two maps are clean if one of them is a submersion.
The main result of this subsection  is the following

\begin{them}
\label{thm:bimodule}
Let $(G\toto G_0, \omega_G +\Omega_G)$, $(H\toto H_0, \omega_H +\Omega_H)$,
  and $(K\toto K_0, \omega_K +\Omega_K )$ be quasi-symplectic groupoids. 
Assume that  $G_{0}\stackrel{\rho_1 }{\leftarrow} X
\stackrel{\sigma_1}{\rightarrow} H_{0}$ is an Hamiltonian
 $G$-$H$-bimodule, and
$H_{0}\stackrel{\rho_2 }{\leftarrow} Y
\stackrel{\sigma_2 }{\rightarrow} K_{0}$
 is an  Hamiltonian $H$-$K$-bimodule.  If 
$Z: = X\times_{H}Y $ is a manifold, then
\begin{enumerate}
\item the two-form $i^* ({\omega_X}\oplus \omega_{Y})
\in \Omega^2 ( X\times_{H_{0}}Y )$, where $i:
X\times_{H_{0}}Y \to X\times Y$ is the  natural inclusion,
descends to a two-form $\omega_Z$ on $Z$; and
\item if moreover assume that 
 $\sigma_1$ and $\rho_2$ are clean, then
   $(Z, \omega_Z )$, equipped with the left $G$-
and right $K$-actions as in Eq. (\ref{eq:GK}),
 is an  Hamiltonian $G$-$K$-\bimodule.
\end{enumerate}
\end{them}
\begin{pf} 
First, note that for any $(x, y)\in  X\times_{H_{0}}Y$,
the tangent space to the $H$-orbit is spanned by vectors
of the form $(\hat{\xi}(x),\hat{\xi}(y) )$, $\forall \xi\in A_H |_m$,
where $A_H$ is the Lie algebroid of
$H$,  and $m=\sigma_1 (x)=\rho_2 (y)$. Here 
we let $H$ act on $X$ from the left: $h\cdot x=xh^{-1}$,
and $\hat{\xi}(x)$ denotes the infinitesimal vector field generated
by this action.
Now
\be
&&(\hat{\xi}(x), \hat{\xi}(y))\per i^* ({\omega_X}\oplus \omega_{Y})\\
&=&\hat{\xi}(x) \per \omega_X  +\hat{\xi}(y) \per \omega_{Y} \\
&=&-\sigma_1^* \epsilon^* (\Vec{\xi} (m) \per \omega_H )
+\rho^*_2 \epsilon^* (\Vec{\xi} (m) \per \omega_H )\\
&=&0.
\ee 

Secondly, let $\lL$ be any local bisection of $H\toto H_0$.
 Then $\lL$  induces a local diffeomorphism on both    $X$ and $Y$, denoted
by $\Phi_\lL$. By the left multiplication, $\lL$ also induces
a local diffeomorphism on $H$ itself, which again, by abuse
of notation,  is  denoted by $\Phi_\lL$.
We need to prove that 
\begin{equation}
\label{eq:inv}
\Phi_\lL^* [ i^* ({\omega_X}\oplus \omega_{Y})]=
i^* ({\omega_X}\oplus \omega_{Y}).
\end{equation}

Given any tangent vectors $(\delta_x^i , \delta_y^i)
\in T_{(x, y)}(X\times_{H_{0}}Y)$,  $i=1, 2$,
let $u^i =\sigma_{1*} \delta_x^i 
=\rho_{2*}\delta_y^i\in T_m H_0$, where $m=\sigma_{1}(x)=\rho_{2}(y)$,
and $\delta_h^i =\Phi_{\lL *}u^i \in T_h H$. 
It is simple to see
that
$(0, \delta_h^i  , \delta_x^i ,  \Phi_{\lL*} \delta_h^i )\in T\Lambda_1$
and $(\delta_h^i ,  0 , \delta_y^i ,  \Phi_{\lL*} \delta_y^i )
\in T\Lambda_2$, where 
$\Lambda_1 \subset G\times \overline{H} \times X\times \overline{X}$
and $\Lambda_2  \subset H\times \overline{K}\times Y\times \overline{Y}$
 are the  corresponding graphs of the groupoid actions.
From the compatibility condition, it follows that
$$ -\omega_H (\delta_h^1 , \delta^2_h)+
\omega_X (\delta_x^1   , \delta^2_x )
-\omega_X (\Phi_{\lL*} \delta_x^1 , \Phi_{\lL*} \delta^2_x ) =0,$$
and
$$ \omega_H (\delta_h^1 , \delta^2_h)+\omega_Y (\delta^1_y , \delta^2_y )
-\omega_Y (\Phi_{\lL*} \delta^1_y , \Phi_{\lL*} \delta^2_y ) =0.$$
Thus we have
$$(\omega_X \oplus \omega_Y ) ( (\delta_x^1, \delta_y^1) , \  (\delta^2_x
 ,\delta^2_y ) )
=(\omega_X \oplus \omega_Y )
((\Phi_{\lL*} \delta_x^1,
\Phi_{\lL*} \delta_y^1 ) ,  \ (\Phi_{\lL*} \delta^2_x , 
\Phi_{\lL*} \delta_y^2 )).$$
Eq. (\ref{eq:inv}) thus follows.  Therefore we conclude that there is a 
two-form $\omega_Z$ on $Z: =X\times_H Y$ such that 
$$\pi^* \omega_Z =i^* ({\omega_X}\oplus \omega_{Y}),$$
where $\pi : X\times_{H_{0}}Y\to Z$ is the projection.

It is straightforward  to check that  
$$d\omega_Z =(\rho_3 \times \sigma_3)^* [\Omega_G\oplus 
\overline{\Omega_K}]$$
and  the two-form $\omega_Z$ is 
 compatible with the  action of the quasi-symplectic
groupoid $G\times \overline{K}\toto G_0 \times \overline{K_0}$.

It remains to prove the minimal  non-degeneracy condition. First we need
the following 

 \begin{lem}
\label{lem:rho}
Let $({G}\toto G_0 , \omega_G +\Omega_G )$
and $({H}\toto H_0 , \omega_H +\Omega_H )$ be quasi-symplectic groupoids,
$G_{0}\stackrel{\rho }{\leftarrow} X
\stackrel{\sigma}{\rightarrow} H_{0}$ an Hamiltonian
$G$-$H$-bimodule  with $\omega_X\in \Omega^2 (X)$.
Then 
\begin{eqnarray}
&&(\ker \rho_* )^{\omega_X}=\{\hat{\xi}(x)|\xi \in A_G|_{\rho (x)} \}+
\ker \omega_X, \ \ \ \mbox{and} \label{eq:23}\\
&&(\ker \sigma_* )^{\omega_X}=\{\hat{\xi}(x)|\xi \in A_H|_{\sigma (x)} \}+
\ker \omega_X , \label{eq:24}
\end{eqnarray}
where $A_G$ and $A_H$  denote the Lie algebroids of $G$ and $H$,
 respectively. 
\end{lem}
\begin{pf}
It is obvious that $\{\hat{\xi}(x)|\xi \in A_G|_{\rho (x)}\}+
\ker \omega_X \subseteq (\ker \rho_* )^{\omega_X}$.
Now 
\be
\dim (\ker \rho_* )^{\omega_X}&=& \dim X-\dim \omega_X^b (\ker \rho_* )\\
&=&\dim \rho_* (T_x X)+\dim \ker \rho_*  -\dim \omega_X^b (\ker \rho_* )\\
&=&\dim \rho_* (T_x X)+ \dim (\ker  \omega_X  \cap \ker \rho_* ) 
\ \ \mbox{(by Lemma \ref{lem:Ga-space} (1))}\\
&\leq&\rank A_G-\dim ( A_G|^{x}_{x}) + \dim (\ker  \omega_X \cap \ker \rho_* )\\
&=&\dim \{\hat{\xi}(x)|\xi \in A_G|_{\rho (x)}\}+\dim (\ker  \omega_X \cap \ker \rho_* ).
\ee

On the other hand, using  Lemma \ref{lem:Ga-space} (2),
 it is easy to check that
\begin{equation}
\label{eq:direct}
(\ker \omega_X \cap \ker  \rho_*) \oplus (\ker \omega_X \cap\{\hat{\xi}(x)|\xi \in A_G|_{\rho (x)}\} )
 =\ker  \omega_X. 
\end{equation}

To prove this equation, first one easily sees that
 $\ker  \omega_X $ can be written as the sum
of the two subspaces on the left hand side.
To show  that this is a direct sum,
 it suffices to show that  the intersection of these two
subspaces is zero. This is because
\be
&&(\ker \omega_X \cap \ker  \rho_* ) \cap  \{\hat{\xi}(x)|\xi \in A_G|_{\rho (x)}\} \\
&\subseteq &\ker \omega_X \cap \ker  \rho_* \cap  \ker\sigma_*\\
&=&\ker \omega_X \cap \ker (\rho\times \sigma )_* \ \ 
 \ \mbox{(by Lemma \ref{lem:Ga-space} (2)b)}\\
&=&0.
\ee

From Eq. (\ref{eq:direct}), it follows that 
\be
&&\dim(\{\hat{\xi}(x)|\xi \in A_G|_{\rho (x)} \}+ \ker \omega_X)\\
&= &\dim \{\hat{\xi}(x)|\xi \in A_G|_{\rho (x)}\} +
\dim \ker \omega_X -\dim( \ker \omega_X \cap \{\hat{\xi}(x)|\xi \in A_G|_{\rho (x)}\})\\
&=&\dim \{\hat{\xi}(x)|\xi \in A_G|_{\rho (x)}\}+\dim (\ker  \omega_X \cap \ker \rho_* ). 
\ee
Thus Eq. (\ref{eq:23}) follows immediately. Similarly
Eq. (\ref{eq:24}) can be proved.
This concludes the proof of the lemma. 
\end{pf}

Assume that $[(\delta_x,   \delta_y) ]\in T_{[(x, y)]}Z$,
where  $(\delta_x,   \delta_y) \in T_{(x, y)}(X\times_{H_{0}}Y )$,
 is in the kernel of $\omega_Z$. Then
\begin{equation}
\label{eq:13}
  \omega_X (\delta_x , \delta'_x )+ \omega_Y ( \delta_y , \delta'_y )=0,
\ \ \ \forall (\delta'_x ,  \delta'_y ) \in   T_{(x, y)}(X\times_{H_{0}}Y).
\end{equation}
By letting $ \delta'_y=0$, it follows that
$\omega_X (\delta_x , \delta'_x )=0$ for any $\delta'_x \in \ker \sigma_{1*}$.
Therefore, according to Lemma \ref{lem:rho}, we have 
$$\delta_x \in (\ker \sigma_{1*})^{\omega_X} =
\{\hat{\eta }(x)|\eta \in A_H|_{\sigma_1 (x)} \}+ \ker \omega_X.           $$
It thus follows that we can always write
$\delta_x =\hat{\xi}(x)+\hat{\eta_1}(x)$
for some $\xi\in A_G|_{\rho_1 (x)}$ and $\eta_1 \in
A_H|_{\sigma_1 (x)}$ such that  $\Vec{\xi}(\rho_1 (x) ) \in \ker \omega_{G}$.

Similarly, one shows that  $\delta_y = \hat{\eta_2} (y)+\hat{\zeta}(y)$,
for some $\eta_2\in A_H|_{\rho_2 (y)}$ and $\zeta \in
A_K|_{\sigma_2 (y)}$ such that $\Vec{\zeta} (\sigma_2 (y)) \in \ker \omega_{K}$.

Now $\sigma_{1*}\delta_x = -a_{A_H} (\eta_1 )$ and $\rho_{2*}\delta_y
=a_{A_H} (\eta_2 )$. Thus we have $\eta_1 -\eta_2 \in \ker a_{A_H}$. From
 Eqs. (\ref{eq:13}) and  (\ref{eq:mom}), it follows that
$$\omega_G (\Vec{\xi}(m) , \rho_{1*}\delta'_x )
-\omega_H (\Vec{\eta_1}(n), \sigma_{1*}\delta'_x )
+\omega_H ( \Vec{\eta_2}(n), \rho_{2*}\delta'_y )
-\omega_K (\Vec{\zeta} (p), \sigma_{2*}\delta'_y )=0, $$
where $m=\rho_{1}(x)$, $n=\sigma_1 (x)=\rho_2 (y)$ and
$p=\sigma_2 (y)$.
Hence $\omega_H  (\Vec{\eta_1-\eta_2 }(n), \delta_n )=0$ for any
$\delta_n\in  f_* T_{(x, y)}(X \times_{H_0} Y )$,
   where $f: X \times_{H_0} Y\to H_0$ is the map 
 $f(x ,y )=\sigma_1 (x)$. By  the clean assumption, we may assume
that $f_* T_{(x, y)}(X \times_{H_0} Y ) =\sigma_{1*}(T_x X)$ (or
$\rho_{2*} (T_{y}Y)$, in which case, a similar proof can be 
carried out). Thus
we have  $\omega_H  (\Vec{\eta_1-\eta_2 }(n), \sigma_{1*}(T_x X))=0$,
which implies that $\hat{\eta_1}(x)-\hat{\eta_2}(x) \in \ker \omega_X^b$.
On the other hand, since 
$(\rho_1 \times \sigma_1 )_* (\hat{\eta_1}(x)-\hat{\eta_2}(x) )
=(0, a_{A_H} (\eta_1 -\eta_2 ))=0$, we have 
$\hat{\eta_1}(x)-\hat{\eta_2}(x)=0 $ according to Lemma
\ref{lem:Ga-space} (2)b.
It thus follows that
$$[(\delta_x , \delta_y )]=[(\hat{\xi}(x)+\hat{\eta_1}(x),
\hat{\eta_2}(y)+\hat{\zeta}(y)]= [(\hat{\xi}(x)+\hat{\eta_2}(x),
\hat{\eta_2}(y)+\hat{\zeta}(y)] 
=[(\hat{\xi}(x),  \hat{\zeta}(y))],$$
which implies the minimal non-degeneracy condition.
This completes the proof.
\end{pf}

\subsection{Reduction}

Theorem \ref{thm:bimodule} has many important  consequences.
As an immediate  consequence, we have the following reduction theorem.

\begin{them}
\label{thm:reduction}
Let  $({\gm}\toto{P}, \omega +\Omega )$ be
  a quasi-symplectic
groupoid,
 and  $(X\stackrel{J}{\to} P, \omega_X )$  an Hamiltonian $\gm$-space.
Assume that $m\in P$ is  a regular value of $J$ and $\gm_m^m$
acts on $J^{-1}(m)$ freely and properly, where
$\gm_m^m$ denotes the isotropy group at $m$.
Then $J^{-1}(m)/\gm_m^m$ is a symplectic manifold.
More generally, if $({\gm_i}\toto{P_i}, \omega_i +\Omega_i )$, $i=1, 2$,
are quasi-symplectic groupoids, 
$(X\stackrel{J_1 \times J_2}{\lon } P_1 \times P_2, \omega_X )$
  is an Hamiltonian $\gm_1 \times \gm_2$-space,  $m\in P_2$ is a regular
value for $J_2 :X\lon P_2$,  and
$(\gm_{2})_m^m$ acts on $J_2^{-1}(m)$ freely and properly,
 then $J_2^{-1}(m)/(\gm_{2})_m^m$ is
  naturally an Hamiltonian $\gm_1$-space.
\end{them}               
\begin{pf}
Note that $(X\stackrel{J_1 \times J_2}{\lon} P_1 \times P_2, \omega_X )$
  being an Hamiltonian $\gm_1 \times \gm_2$-space is equivalent to
 $X$ being  a $\gm_1$-$\overline{\gm_2}$-bimodule by considering
$X$ as  a right $\gm_2$-space.
Let $\O \subset P_2 $ be the groupoid orbit of $\gm_2$
 through $m$. Then
$P_2\stackrel{i}{\leftarrow}\bar{\O}
\stackrel{}{\rightarrow}\cdot $
is  an Hamiltonian $\overline{\gm_2}$-$\cdot$-bimodule
according to Proposition \ref{prop:example}.
The clean assumption is satisfied since $J_2: J_2^{-1}(\O)
\to \O$ is a submersion.
By Theorem \ref{thm:bimodule}, $X\times_{\overline{\gm_2}}\bar{\O}$
is an Hamiltonian $\gm_1$-$\cdot$-Hamiltonian bimodule, i.e., a
Hamiltonian $\gm_1$-space. It is easy to see
 that $X\times_{\overline{\gm_2}}\bar{\O}$
is naturally diffeomorphic to  $J_2^{-1}(m)/(\gm_{2})_m^m$.
\end{pf}

\begin{numrmk}
As a consequence, $X/\gm$ (assuming  being a smooth
manifold) is naturally a Poisson manifold.
One should also be able to see this using the
reduction of Dirac structures, as an Hamiltonian
$\gm$-space infinitesimally corresponds to
some particular Dirac structure \cite{BCWZ}.
\end{numrmk} 

Various reduction theorems in the literature are indeed consequences
of Theorem \ref{thm:reduction}. In particular, 
applying Theorem \ref{thm:reduction}  to the AMM quasi-symplectic
groupoids, we recover the Hamiltonian reduction  theorem
of quasi-Hamiltonian $G$-spaces of Alekseev--Malkin--Meinrenken
\cite{AMM}.

\begin{cor}
Let $X$ be a quasi-Hamiltonian $G_1\times G_2$-space and let
$f\in G_2$ be a regular value of the momentum
 map $J_2: X\rightarrow G_2$.
Then the pull-back of the 2-form $\omega$ to $J_2^{-1}(f)$
descends to the {\em reduced space}
\begin{equation*}
X_f = J_2^{-1}(f)/(G_2)_f
\end{equation*}
(assuming it is a smooth manifold) and makes it into a quasi-Hamiltonian
 $G_1$-space. Here $(G_2)_f$ is the 
isotropy group of $G_2$ at $f$.
 In particular, if $G_1=\{e\}$ is trivial then $X_f$ is a symplectic manifold. 
\end{cor}

Another immediate  consequence of  Theorem \ref{thm:bimodule} is the 
following

\begin{them}
Let   $({\gm}\toto{P}, \omega +\Omega )$ be  a quasi-symplectic
groupoid, and  $(X\stackrel{J_1}{\to} P, \omega_X )$,
and $(Y\stackrel{J_2}{\to} P, \omega_Y )$ be  Hamiltonian $\gm$-spaces.
Assume that $J_1 :X\to P$ and $J_2 : Y\to P$ are clean.
Then  $X\times_\gm \overline{Y}$ is a symplectic manifold.
\end{them}
We will call $X\times_\gm \overline{Y}$ the classical intertwiner
 space between $X$ and $Y$.
When ${\gm}\toto{P}$ is the  symplectic groupoid $T^*G\toto \frakg^*$,
this reduces to  the  classical intertwiner space
$(X\times \overline{Y})_0$ of Hamiltonian $G$-spaces \cite{GS}.
We refer the reader to \cite{Xu:intert}
for the detailed study of classical intertwiner spaces
of  symplectic groupoids.

\section{Morita equivalence}

This section is devoted to the study of Morita equivalence
of quasi-symplectic groupoids. The main result is
that Morita equivalent quasi-symplectic groupoids define
equivalent momentum map theories. See 
Theorem \ref{thm:bijection}     and Corollary \ref{cor:reduction}.

\subsection{Morita equivalence of quasi-symplectic groupoids}

Morita equivalence is an important equivalence relation
for groupoids. Indeed  groupoids moduli Morita equivalence
can be identified with the so called stacks, which are
useful in the study of singular spaces such as moduli spaces.
Morita equivalence of symplectic groupoids were studied in
\cite{Xu:90}. Here we will generalize this notion to
quasi-symplectic groupoids. Let us first recall the
definition of Morita equivalence of Lie groupoids \cite{LTX, Xu:90}.

\begin{defn}
 Lie groupoids $G\rightrightarrows G_0$ and
$H\rightrightarrows H_{0}$ are said to be
 Morita equivalent if  there exists a
manifold $X$ together with two surjective submersions
$$G_{0}\stackrel{\rho }{\leftarrow} X
\stackrel{\sigma}{\rightarrow} H_{0},$$
a left action of $G$ with respect to $\rho$,
a right action of $H$ with respect to $\sigma$ such that
\begin{enumerate}
 \item the two actions commute with each other;
 \item $X$ is a locally trivial $G$-principal
bundle over $ X \stackrel{\sigma}{\rightarrow} H_{0}$; and
\item  $X$ is a locally trivial $H$-principal
bundle over  $G_{0}\stackrel{\rho}{\leftarrow} X$.
\end{enumerate}  
In this case,  
$G_0 \stackrel{\rho}{\leftarrow} X\stackrel{\sigma}{\to} H_0 $
is called an  equivalence bimodule between the  Lie groupoids $G$ and $H$.
\end{defn}

It is known that de-Rham cohomology groups are invariant
under  Morita equivalence.
 
\begin{prop}
 \cite{BEFFGK, BX,  hae84}
If  $G\toto G_{0}$ and $H\toto H_{0}$
are  Morita equivalent Lie groupoids, then
$$ H_{DR}^k(G\lcom) \longiso H_{DR}^k({H}\lcom)$$
\end{prop}

\begin{defn}
Quasi-symplectic groupoids $({G}\toto G_0 , \omega_G +\Omega_G )$
and $({H}\toto H_0 , \omega_H +\Omega_H )$ are  said to
be Morita equivalent
if there exists  a  Morita equivalence bimodule
 $G_{0}\stackrel{\rho }{\leftarrow} X
\stackrel{\sigma}{\rightarrow} H_{0}$
between the Lie groupoids $G$ and $H$,
together with a two-form $\omega_X\in \Omega^2 (X)$ such that
$X$ is also an  Hamiltonian $G$-$H$-\bimodule.
\end{defn}

Suppose that $G \toto G_0$ and $H \toto H_0 $
 are Morita  equivalent Lie groupoids with equivalence bimodule
$G_0 \stackrel{\rho}{\leftarrow} X\stackrel{\sigma}{\to} H_0$.
  We say that $m\in G_{0}$ and $n\in H_{0}$ are
related if $\rho^{-1}(m)\cap \sigma^{-1}(n)\neq \emptyset $. 
The following are basic properties \cite{Xu:90}.

\begin{prop}
If $G \toto G_0$
and $H \toto  H_0$ are Morita equivalent 
Lie groupoids with equivalence bimodule 
$G_0 \stackrel{\rho}{\leftarrow} X\stackrel{\sigma}{\to} H_0$.
Assume that $m_0\in G_0$ and $n_0\in H_0$ are
related.  Then
\begin{enumerate}
 \item $\dim G + \dim H = 2\dim X$;
\item an  element $n\in H_{0}$ is related to $m_{0}$ if and only if $n$ lies in the same
groupoid orbit as $n_{0}$; and conversely, $m\in G_{0}$ is 
related to $n_{0}$ if and only if $m$ lies in the same groupoid orbit
 as  $m_{0}$; and
 \item the isotropy groups at $m_{0}$ and $n_{0}$ are isomorphic.
\end{enumerate}           
\end{prop}

\begin{them}
Morita equivalence is indeed  an  equivalence relation for 
quasi-symplectic groupoids.
\end{them}
\begin{pf}
From Proposition \ref{prop:example} (1), we know that Morita equivalence is reflective.
If $G_{0}\stackrel{\rho }{\leftarrow} X
\stackrel{\sigma}{\rightarrow} H_{0}$   is an Hamiltonian bimodule
defining the Morita equivalence between $(G\toto G_0, \omega_G +\Omega_G)$
and $(H\toto H_0, \omega_H +\Omega_H)$,  then
$H_{0}\stackrel{\sigma }{\leftarrow} \overline{X} \stackrel{\rho}{\rightarrow} 
G_{0}$, 
with the reversed actions,  is an Hamiltonian bimodule
defining the Morita equivalence between $(H\toto H_0, \omega_H +\Omega_H)$
and $(G\toto G_0, \omega_G +\Omega_G)$.
So the symmetry  follows.
As for the transitivity, let
 $(G\toto G_0, \omega_G +\Omega_G)$, $(H\toto H_0, \omega_H +\Omega_H)$,
  and $(K\toto K_0, \omega_K +\Omega_K )$
be quasi-symplectic groupoids.  Assume that
  $G_{0}\stackrel{\rho_1 }{\leftarrow} X
\stackrel{\sigma_1}{\rightarrow} H_{0}$ is a $G$-$H$ equivalence 
bimodule, and $H_{0}\stackrel{\rho_2 }{\leftarrow} Y
\stackrel{\sigma_2 }{\rightarrow} K_{0}$ an   $H$-$K$-equivalence
 bimodule,  respectively. 
It is known that $Z=X\times_H Y$ is a bimodule
defining the Morita equivalence between the groupoids
$G\toto G_0$ and $K\toto K_0$. According to 
Theorem \ref{thm:bimodule}, $Z$ is
also an Hamiltonian $G$-$K$-bimodule. Thus
 $(G\toto G_0, \omega_G +\Omega_G)$
and $(K\toto K_0, \omega_K +\Omega_K )$  are Morita equivalent 
quasi-symplectic groupoids.
\end{pf}

In what follows, we describe some useful constructions
of producing Morita equivalent quasi-symplectic groupoids.

Let $\gm \toto P$ be a  Lie groupoid, and $\omega_i +\Omega_i\in
\Omega^2 (\gm )\oplus \Omega^3 (P )$, $i=1, 2$,  be two cohomologous
3-cocycles. This means that
 there are  $B\in \Omega^2 (P )$ and $\theta \in \Omega^1 (\gm )$
such that 
$$(\omega_1+\Omega_1)-(\omega_2+\Omega_2) =\delta (B+\theta ). $$
Following \cite{BR}, we say that $\omega_1+\Omega_1$ and 
$\omega_2+\Omega_2$ differ by a  {\em gauge transformation of the first 
type} if $(\omega_1+\Omega_1)-(\omega_2+\Omega_2) =\delta B$, i.e.,
$$\omega_1-\omega_2 =s^* B-t^* B, \ \ \ \Omega_1 -\Omega_2 =d B. $$
And we say that   $\omega_1+\Omega_1$ and
$\omega_2+\Omega_2$ differ by a  {\em gauge transformation of the second 
type} if $(\omega_1+\Omega_1)-(\omega_2+\Omega_2) =\delta \theta $, i.e., 
$$ \Omega_1=\Omega_2,
\ \ \ \omega_1 =\omega_2 -d\theta, \ \ \ \mbox{for some }\theta \in \Omega^1 (\gm)
\ \ \mbox{such that } \partial \theta =0.$$ 

It is simple to see that gauge transformations of the first type transform
quasi-symplectic groupoids into quasi-symplectic groupoids (see also
 \cite{BCWZ}).
Below we see that the resulting quasi-symplectic groupoids are indeed
Morita equivalent (see  \cite{BR} for the case
of symplectic groupoids).

\begin{prop}
\label{prop:4.7}
Assume that $(\gm\toto P, \omega +\Omega )$ is
 a quasi-symplectic groupoid.
Then $(\gm\toto P, \omega' +\Omega' )$, where $\omega' =\omega +s^* B
-t^*B$ and $\Omega' =\Omega +dB$, for any $B\in \Omega^2 (P)$,  is a
Morita equivalent quasi-symplectic groupoid.
\end{prop}
\begin{pf}
First, we need to show that $\omega'$ is non-degenerate in the sense
of Definition \ref{def:non-deg}. By Proposition 
\ref{prop:eq}, it suffices to show that
$a: \ker \omega'_m \cap  {A}_{m}\to \ker \omega'_m \cap T_{m}P$ is
injective. Assume that $\xi \in \ker \omega'_m \cap  {A}_{m}$ such that
$a(\xi )=0$.  Then we have for any $v\in T_m P$, 
 $$0=\omega' (\Vec{\xi}, v )=(\omega +s^*B -t^*B)  (\Vec{\xi}, v)
=\omega (\Vec{\xi}, v)+B(a(\xi), v)=\omega (\Vec{\xi}, v). $$  
Thus we have  $\xi \in \ker \omega_m \cap  {A}_{m}$, which implies that
$\xi =0$.

To prove the Morita equivalence, let $X=\gm$ and $\omega_X=
\omega +s^*B$. We let  $(\gm\toto P, \omega' +\Omega' )$
act on $X$ from the left by left multiplications
 and let $(\gm\toto P, \omega +\Omega )$
act on $X$ from the right by right multiplications.
 It is simple to check that these actions
are compatible with  the quasi-symplectic structures.
 It remains to check the minimal  non-degeneracy
condition.  Assume that $\delta_x \in \ker \omega_X$. Then
for any $\zeta \in A_{t(x)}$, we have,  
$$0=\omega_X (\delta_x, \ceV{\zeta } (x))=\omega (\delta_x, \ceV{\zeta } (x))
+B(s_* \delta_x, s_* \ceV{\zeta })=\omega (\delta_x, \ceV{\zeta } (x)) $$
 since $s_* \ceV{\zeta } (x)=0$.
Hence
$\omega (t_* \delta_x,  \ceV{\zeta }(t(x)))=0$ according to 
Eq. (\ref{eq:left-right}), which   implies that
$t_* \delta_x\in \ker \omega $. Therefore $t_* \delta_x =a (\eta )$
for some $\eta \in A_{t(x)}$ such that
 $\Vec{\eta } (t(x)) \in \ker \omega$. Hence
 $\ceV{\eta } (t(x)) \in \ker \omega$ by Corollary \ref{cor:2.3} (3),
which implies that  $\ceV{\eta } (x) \in \ker \omega$.
Set $\delta'_x =\delta_x +\ceV{\eta} (x)$. Thus
$t_* \delta'_x =t_* \delta_x +t_* \ceV{\eta} (x)
=t_* \delta_x -a(\eta )=0$. Hence we may write
$\delta'_x =\Vec{\xi}(x)$ for some $\xi
\in A_{s(x)}$. Moreover, a simple computation yields that 
$$\Vec{\xi}(x)\per \omega'= \delta'_x \per \omega'
=\delta'_x \per \omega_X -\delta'_x \per t^* B=\delta'_x \per \omega_X
=\delta_x \per \omega_X+\ceV{\eta} (x)\per  \omega+ \ceV{\eta} (x)\per s^*B
=0.$$
Thus $\Vec{\xi}(s(x))\in \ker \omega'$ according to 
Eq. (\ref{eq:left-right}). This concludes the proof.  
\end{pf}

\begin{numrmk}
Note that quasi-symplectic groupoids are  in general  not preserved
under gauge transformations of the second type. For instance,
the symplectic structure $\omega$ on the symplectic
groupoid $T^* G\toto \frakg^*$ is $d\theta$, where $\theta\in \Omega^1
(T^* G)$ is the Liouville one-form. It is simple to see that
 $\theta$ satisfies
the condition  $\partial \theta =0$. However $T^* G\toto \frakg^*$ with 
the zero two-form is clearly not quasi-symplectic.
\end{numrmk}

For a Lie groupoid $\gm \toto P$ and  a surjective
submersion $\phi: Y\to P$, 
 we denote by $\Gamma[Y]$ the subgroupoid of $(Y\times Y)\times \Gamma$
consisting of 
$\{(y_1,y_2, r)|\; s(r) =\phi (y_1 ), \ t(r)=\phi (y_2 ) \}$,
called the pull-back groupoid.
Clearly the
 projection $\pr: \Gamma[Y]\to \gm$ defines  a groupoid homomorphism.
By abuse   of  notations, we also use $\pr$ to denote
the corresponding map on the unit spaces $\phi: Y\to P$.

\begin{prop}
\label{prop:4.9}
Assume that $(\gm\toto P, \omega +\Omega )$ is
 a quasi-symplectic groupoid, and $\phi: Y\to P$ a surjective
submersion.
Then $(\gm [Y]\toto Y, \pr^*\omega + \pr^* \Omega )$ is a 
quasi-symplectic groupoid Morita equivalent to $(\gm\toto P, \omega +\Omega )$.
\end{prop}
\begin{pf}
It is obvious  that  $\pr^*\omega + \pr^* \Omega$ is a 3-cocycle
since $\pr$ is a Lie groupoid homomorphism. By $A_Y$, we denote the
Lie algebroid of $\gm [Y]\toto Y$. It is simple to see that $\forall y
\in Y$,
$$A_Y|_y =\{(\delta_y, \xi )|\delta_y \in T_y Y, \ \xi\in A_{\phi (y)}
\ \mbox{ such that }\  \phi_* \delta_y= a(\xi )\},$$
 with the anchor $a_Y: A_Y\to TY$ being  given 
by the projection $(\delta_y , \xi )\to \delta_y$,
 where $A$ is the Lie algebroid of $\gm$.  Therefore,  an element
$(\delta_y , \xi ) \in A_Y|_y$, where $\phi_* \delta_y=
a(\xi )$, belongs to  $  \ker (\pr^* \omega )|_y \cap  A_Y|_y$ if and only if
 $ \xi\in \ker \omega_{\phi (y)}\cap A_{\phi (y)}$.
 This  implies
that $a_Y: \ker (\pr^* \omega ) \cap A_Y \to 
\ker (  \pr^* \omega ) \cap TY$ is indeed injective.
It thus follows that $ (\gm [Y]\toto Y, \pr^*\omega + \pr^* \Omega )$
is a quasi-symplectic groupoid by  dimension counting.

To show the Morita equivalence, let $X:= \gm \times_{t, P, \phi} Y$
and $\omega_X =p^* \omega$, where $p: X\to \gm$ is the
natural projection. It is standard  that 
$P\stackrel{\rho }{\leftarrow} X \stackrel{\sigma }{\rightarrow} Y$
is a  $\gm$-$\gm [Y]$-bimodule defining a  Morita equivalence
between these two Lie  groupoids, where 
$$\rho ( r, y)=s (r) , \  \ \mbox{ and } \sigma (r, y)=y$$
and the left $\gm$-action is
\begin{equation}
\label{eq:action1}
\tilde{r} \cdot (r, y)=(\tilde{r} r, y), 
\ \ \ \ t(\tilde{r} )=s (r)
\end{equation}
while the right $\gm [Y]$-action is

\begin{equation}
\label{eq:action2} 
(r, y)\cdot (y_1, y_2, \tilde{r} )=(r\tilde{r}, y_2 ), \ \ \ y=y_1, \ \ t(r)=\phi (y)=\phi (y_1 )=s (\tilde{r} ).
\end{equation}

 It is also simple to
check that $\omega_X$ is compatible with 
the $\gm\times \overline{\gm [Y]}$-action. For the 
minimal  non-degeneracy condition, assume that
$(\delta_r , \delta_y )\in T_{(r, y)}X$ such that
$(\delta_r , \delta_y )\per \omega_X=0$, which is
equivalent to that $\delta_r  \per \omega =0$.
By Proposition \ref{prop:example},  we have 
$\delta_r =\Vec{\xi}(r)-\ceV{\eta}(r)$,
where $\xi \in A_{s(r)}$ and $\eta 
\in A_{t(r)}$ such that 
 $ \Vec{\xi}( s(r))$ and  $\Vec{\eta} (t (r) )\in
\ker \omega $. Thus
$(\delta_r , \delta_y ) =\hat{\xi}(r, y)-\hat{\eta'}(r, y)$,
where  $\eta' =(\delta_y , \eta )\in A_{Y}|_{y}$
 clearly satisfies the
condition that $\Vec{\eta'} (y )\in \ker  \pr^*\omega$.
This concludes the proof.
\end{pf}

Combining  Proposition \ref{prop:4.7} and Proposition  \ref{prop:4.9},
 we are lead to the following

\begin{them}
\label{thm:4.15}
Let  $({G}\toto G_0 , \omega_G +\Omega_G )$
and $({H}\toto H_0 , \omega_H +\Omega_H )$ be
pre-quasi-symplectic groupoids, which are Morita
 equivalent as Lie groupoids with
an equivalence bimodule $G_{0}\stackrel{\rho }{\leftarrow} X
\stackrel{\sigma}{\rightarrow} H_{0}$.
If $\rho^*  (\omega_G +\Omega_G )$ and
 $\sigma^* (\omega_H +\Omega_H )$,
as de-Rham  3-cocycles of the groupoid $G[X]\cong H[X]\toto X$,
differ by a gauge transformation of the first type, then
if one is quasi-symplectic, so is the other. Moreover, they
are Morita equivalent as quasi-symplectic groupoids.
\end{them}

\subsection{Generalized homomorphisms of quasi-symplectic
groupoids}

Recall  that   a {\em generalized homomorphism}
from  a  Lie groupoid $G\toto G_0$ to $H\toto  H_{0}$ 
 is given by a manifold $X$, two smooth maps
$G_{0}\stackrel{\rho }{\leftarrow} X
\stackrel{\sigma}{\rightarrow} H_{0}$,
a left action of $G$ with respect to $\rho$,
a right action of $H$ with respect to $\sigma$, such that the two
actions commute, and   
$X$ is a locally trivial $H$-principal
bundle over $G_{0}\stackrel{\rho}{\leftarrow}X$ 
\cite{LTX}. In particular, $\rho:  X\to G_{0}$ 
must be a surjective submersion, and the 
(right) $H$-action on $X$ is free and proper.

Generalized homomorphisms can be composed  just like the usual
groupoid homomorphisms; 
 thus there is a
category ${\mathcal{G}}$ whose objects are Lie groupoids and morphisms
are generalized homomorphisms \cite{hae84,hilsum-skandalis87,TXL},
where isomorphisms in the category ${\mathcal{G}}$ are just Morita
equivalences \cite{muhly-renault-williams87, Xu:90}.

Similarly, we can introduce the notion of  generalized homomorphisms
between quasi-symplectic groupoids.

\begin{defn}
 A  generalized homomorphism from a quasi-symplectic groupoid
$({G}\toto G_0 , \omega_G +\Omega_G )$ to 
a quasi-symplectic groupoid  $({H}\toto H_0 , \omega_H +\Omega_H )$
is  an  Hamiltonian $G$-$H$-bimodule
$G_{0}\stackrel{\rho }{\leftarrow} X
\stackrel{\sigma}{\rightarrow} H_{0}$,  which is, in the same
time,  also a generalized homomorphism  from $G$ to $H$.
\end{defn}

Theorem \ref{thm:bimodule} implies the following:

\begin{them}
\label{thm:4.7}
There is a category,   whose objects are quasi-symplectic
groupoids, and morphisms
are generalized homomorphisms of quasi-symplectic
groupoids.
The isomorphisms in this  category correspond exactly to  Morita
equivalences  of quasi-symplectic groupoids.
\end{them}

It is known that a strict homomorphism of Lie groupoids must be 
a generalized homomorphism.
For quasi-symplectic groupoids, one can also introduce 
the notion of strict homomorphisms.

\begin{defn}
\label{defn:4.8}
A strict homomorphism of quasi-symplectic groupoids
from  $({G}\toto G_0 , \omega_G +\Omega_G )$  to
$({H}\toto H_0 , \omega_H +\Omega_H )$ is a  Lie
groupoid  homomorphism $\phi: G\to H$ satisfying
\begin{enumerate}
\item $\phi^* (\omega_H +\Omega_H)=\omega_G +\Omega_G$, and 
\item if $\xi \in A_H|_m$ satisfies the properties that
$a_{H}(\xi )=0$ and $\phi^* ( \Vec{\xi} (m) \per \omega_H )=0$,
then $\xi=0$, where $A_H$ is the Lie algebroid of $H\toto H_0$ and
$a_H : A_H\to TH_0$ denotes its anchor map. 
\end{enumerate}
\end{defn}

\begin{prop}
\label{pro:4.9}
For  quasi-symplectic groupoids, strict homomorphisms imply
generalized homomorphisms.
\end{prop}                                            
\begin{pf}
Assume that $\phi: G\to H$ is a   strict homomorphism of quasi-symplectic
groupoids from  $({G}\toto G_0 , \omega_G +\Omega_G )$  to
 $({H}\toto H_0 , \omega_H +\Omega_H )$.
Let $X= G_0 \times_{\phi ,H_0, s} H $,
 and set  $\rho (g_0 ,h)= g_0  $, 
$\sigma(g_0 ,h) = t(h)$. Define a left 
$G$- and a right $H$-action on $X$, respectively, by
 $$g \cdot (g_0 ,h)= ( s (g), \phi (g)h), \ \ \  \mbox{and } (g_0 ,h)\cdot h'
= (g_0 ,  h h').$$
 One checks that this defines  a generalized homomorphism from $G \toto G_0$ 
to $H \toto H_0$ \cite{LTX}. Let $\omega_X =i^* (0, \omega_H)$,
where $i: G_0 \times_{\phi ,H_0, s} H \subset G_0 \times H$ is the
inclusion. It is simple to   see that    $\omega_X $ is compatible with
the $G$-$H$-bi-actions.
It remains to prove the minimal non-degeneracy condition.
Note that for any $\xi \in \gm (A_G )$ and $\eta \in \gm (A_H )$, the vector 
field on $X$ generated by the  infinitesimal action of
$(\xi, \eta)$ is given by 
$$\Hat{(\xi \oplus \eta )} (g_0, h)=
(a_{G} (\xi )(g_0 ), \Vec{\phi_* \xi}(h)-\ceV{\eta}(h) ), \ \ \forall 
(g_0, h)\in X.$$

Assume that $\delta_x =(\delta_{g_0}, \delta_h )\in 
\ker \omega_X$, where $x=(g_0, h)\in X$. This implies
that
\begin{eqnarray}
&&\phi_* \delta_{g_0} =s_* \delta_h \ \ \ (\mbox{thus } 
\phi (g_0 )=s (h)), \ \ \mbox{ and} \\
&&\omega_H (\delta_h, \delta'_h )=0, \ \ \forall \delta'_h 
\in T_h H \ \ \mbox{such that } \ s_* \delta'_h 
\in \Im( \phi_* ). \label{eq:28}
\end{eqnarray}
In particular, for any $\zeta \in A_H |_{t(h)}$, since
$s_* \ceV{\zeta}(h)=0$, which is always in the image of
$\phi_*$, we have  $\omega_H (\delta_h, \ceV{\zeta}(h))=0$. From
 Eq. (\ref{eq:left-right}), it thus follows  that 
$\omega_H (t_* \delta_h , \ceV{\zeta}(t(h)))=0$, which implies that
$t_* \delta_h\in \ker \omega_H$.
By the non-degeneracy assumption of Definition
\ref{def:non-deg}, we have $t_{*}\delta_h =a_H (\eta )$,
 where  $\eta \in A_H|_{t(h)}$  such that
$\Vec{\eta } (t(h)) \in \ker \omega_H $.
 Hence $\ceV{\eta } (t(h)) \in \ker \omega_H $
 according to Corollary \ref{cor:2.3} (3), which
implies that
 $\ceV{\eta}(h)$ belongs to $\ker \omega_H$
 by  Eq. (\ref{eq:left-right}).
 Let $\tilde{\delta}_h =\delta_h +\ceV{\eta}(h)$. Then we have
$t_* \tilde{\delta}_h =t_* \delta_h +t_* \ceV{\eta}(h)=t_* \delta_h 
-a_H (\eta )=0$. Thus $\tilde{\delta}_h   =\Vec{\xi_1} (h)$ for
some $\xi_1\in A_H |_{s(h)}$,
 and hence we have $\delta_h =
\Vec{\xi_1} (h)-\ceV{\eta}(h)$.
On the other hand, for any 
 $\chi \in  A_G|_{g_0}$, since
 $s_* \Vec{\phi_* \chi} (h)
=a_H (\phi_* \chi )=\phi_* a_{G} (\chi )\in \Im \phi_*$, we have
$\omega_H (\delta_h , \Vec{\phi_* \chi} (h) )=0$ by
 Eq. (\ref{eq:28}). Now
 $$\omega_H (\delta_h , \Vec{\phi_* \chi} (h) )=
\omega_H (\Vec{\xi_1} (h)-\ceV{\eta}(h), \Vec{\phi_* \chi} (h) )
=\omega_H (\Vec{\xi_1} (h), \Vec{\phi_* \chi} (h) )
=\omega_H (s_* \Vec{\xi_1} (h), \Vec{\phi_* \chi} (s(h))), $$
where  we used Eq. (\ref{eq:left-right}) in the last equality.

Now $s_* \Vec{\xi_1} (h)= s_* (\delta_h +\ceV{\eta}(h))=s_* \delta_h 
=\phi_* \delta_{g_0}$.
Therefore we have $\omega_H (\phi_* \delta_{g_0}, 
\Vec{\phi_* \chi} (s(h)))=0,
\ \forall \chi \in A_G|_{g_0}$.
It thus follows that $\delta_{g_0}\in \ker (\phi^* \omega_H )$. Since 
$(G\toto G_0, \phi^* \omega_H +\phi^*\Omega_H )$ is 
quasi-symplectic, by the   non-degeneracy assumption, we 
conclude  that $\delta_{g_0}=a_G (\xi)$ for
some $\xi \in A_G|_{g_{0}}$ such that
$\Vec{\xi}(g_0 ) \in \ker (\phi^*\omega_H ) $.
Therefore  for any $\delta'_{g_0}
\in T_{g_0}G_0, \ \omega_H (\Vec{\phi_*\xi} (s(h)),
 \  \phi_* \delta'_{g_0} ) 
=(\phi^* \omega_H )(\Vec{\xi}(g_0), \delta'_{g_0})=0$.
Hence  $\omega_H (\Vec{\xi_1} (s(h))-\Vec{\phi_*\xi} (s(h)), \  \phi_* \delta'_{g_0} )
=\omega_H (\Vec{\xi_1} (s(h)), \  \phi_* \delta'_{g_0} )$.
Since $s: H\to H_0$ is a submersion, we may assume that
$\phi_* \delta'_{g_0}=s_* \delta''_h$ for some  $ \delta''_h\in
T_h H$, and therefore
$$\omega_H (\Vec{\xi_1} (s(h)), \  \phi_* \delta'_{g_0} ) 
=\omega_H (\Vec{\xi_1} (s(h)), \ s_* \delta''_h)
=\omega_H (\Vec{\xi_1}(h), \  \delta''_h)  
=\omega_H ( \delta_h +\ceV{\eta}(h), \ \delta''_h)        
=\omega_H (  \delta_h , \ \delta''_h)   =0$$
by Eq. (\ref{eq:28}) since $ s_* \delta''_h =\phi_* \delta'_{g_0}\in
\Im \phi_*$.
On the other hand, since 
 $$a_{H}(\xi_1 -\phi_*\xi )
=s_{*}\Vec{\xi_1} (h)-a_{H}(\phi_*\xi )=s_{*}\delta_h -\phi_* (a_{G}(\xi) )
=s_{*}\delta_h-\phi_*\delta_{g_0}=0, $$
  it  follows that $\xi_1 -\phi_*\xi=0$. 
Therefore, we conclude that 
$\delta_x =(\delta_{g_0}, \delta_h )=
\Hat{(\xi \oplus \eta )} (g_0, h)$, where  $\Vec{\xi}(g_0 ) \in 
\ker (\phi^*\omega_H ) $ and $ \Vec{\eta}(t(h)) \in \ker  \omega_H $.
\end{pf}

\begin{numrmk}
Note that the second condition in Definition \ref{defn:4.8} is
necessary for Proposition \ref{pro:4.9} to hold.
For instance, given a quasi-symplectic
groupoid $H\toto H_0$ and a fixed point in $H_0$,
one may always think of  this point
 as  a groupoid homomorphism
from $\cdot \toto \cdot $ to $H\toto H_0$.
The first condition is satisfied automatically.
However, $\cdot \stackrel{\rho }{\leftarrow} H
\stackrel{\sigma}{\rightarrow} H_{0}$ is, in general, not
a generalized  homomorphism of quasi-symplectic groupoids
since $H$ is not, in general, an
 Hamiltonian $H$-space under the right $H$-action.
\end{numrmk}

The following proposition describes the
precise relation between generalized homomorphisms
and strict homomorphisms for quasi-symplectic groupoids.

\begin{prop}
Any generalized homomorphism  of quasi-symplectic groupoids
 is equivalent to the composition of a Morita 
equivalence 
 with a strict homomorphism. 
\end{prop}
\begin{pf}
The inverse direction follows from Proposition \ref{pro:4.9}
and Theorem \ref{thm:4.7}, so it remains
to prove the other direction.

Assume that $G_{0}\stackrel{\rho }{\leftarrow} X
\stackrel{\sigma}{\rightarrow} H_{0}$
is a generalized homomorphism of
quasi-symplectic groupoids from
$({G}\toto G_0 , \omega_G +\Omega_G )$   to
$({H}\toto H_0 , \omega_H +\Omega_H )$.

Consider the transformation groupoid
$Q:=(G\times H)\times_{(G_0\times H_0 )} X\toto X $
as in Proposition \ref{pro:pull}.
One easily checks that $Q\toto X$ is
isomorphic to $G[X]\toto X$, where the isomorphism
is given by $(g, h, x) \to (x, g^{-1}xh , g)$, $\forall
(g, h, x)\in (G\times H)\times_{(G_0\times H_0)} X$.
Therefore we have two groupoid
homomorphisms $\pr_1 :G[X] \to G$ and
$\pr_2 :G[X] \to H$. 
Equip $G[X]\toto X$ with the 3-cocycle
$$\omega_{G[X]} +\Omega_{G[X]}:= \pr_1^* (\omega_G +\Omega_G)-\delta \omega_X.$$
By Theorem \ref{thm:4.15},
 we know that $(G[X]\toto X, \omega_{G[X]} +\Omega_{G[X]})$
is Morita equivalent to $({G}\toto G_0 , \omega_G +\Omega_G )$.
On the other hand, according to Proposition \ref{pro:pull},
 we have $\omega_{G[X]} +\Omega_{G[X]}=\pr_2^* (\omega_H +\Omega_H )$.
It thus follows from  Theorem \ref{thm:bimodule} that
$X\stackrel{\rho }{\leftarrow}
X\times_{\sigma , H_0, s} H \stackrel{\sigma}{\rightarrow} H_{0} $
is  an Hamiltonian $G [X]$-$H$ bimodule defining 
a generalized homomorphism of quasi-symplectic
groupoids from $(G[X]\toto X, \omega_{G[X]} +\Omega_{G[X]})$ 
to  $({H}\toto H_0 , \omega_H +\Omega_H )$. Here the
two-form $\omega_Z$ on $Z: = X\times_{\sigma , H_0, s} H $ is 
given by $\omega_Z =i^* (0, \omega_H )$, where $i: Z\to X\times
H$ is the inclusion.
By Lemma \ref{lem:Ga-space} 2(b), one easily sees that
Condition (2) in Definition  \ref{defn:4.8} 
is satisfied so that
 $\pr_2 :G[X] \to H$ is  indeed a strict homomorphism of
quasi-symplectic groupoids. This completes the proof.
\end{pf}

The proof of the above proposition also yields the following

\begin{cor}
If  $f: G_{0}\stackrel{\rho }{\leftarrow} X
\stackrel{\sigma}{\rightarrow} H_{0}$ 
is  a  generalized homomorphism of 
quasi-symplectic groupoids 
from $({G}\toto G_0 , \omega_G +\Omega_G )$ to
$({H}\toto H_0 , \omega_H +\Omega_H )$,
then $f^* [\omega_H +\Omega_H ]=[\omega_G +\Omega_G]$,
where $f^*: H^3_{DR} (H\lcom) \to H^3_{DR} (G\lcom)$ is the
induced homomorphism of the de Rham cohomology groups.

In particular, if  $({G}\toto G_0 , \omega_G +\Omega_G )$  and
$({H}\toto H_0 , \omega_H +\Omega_H )$ are Morita equivalent
quasi-symplectic groupoids, then $[ \omega_G +\Omega_G] $ and
$[\omega_H +\Omega_H ] $ define the same class
under the isomorphism $H^3_{DR} (G\lcom) \simeq  H^3_{DR} (H\lcom)$.
\end{cor}

\subsection{Hamiltonian spaces for Morita equivalent
 quasi-symplectic groupoids}
\begin{defn}
Assume that $({G}\toto G_0 , \omega_G +\Omega_G )$
and $({H}\toto H_0 , \omega_H +\Omega_H )$ are Morita equivalent
quasi-symplectic groupoids with
  an equivalence bimodule $G_{0}\stackrel{\rho }{\leftarrow} X
\stackrel{\sigma}{\rightarrow} H_{0}$.
Let $\phi: F\to G_{0}$ be an Hamiltonian $G$-space,
 and $\psi: E\to H_{0}$ an Hamiltonian
$H$-space.  We say that  $F$ and $E$ are a {\sl pair of related
Hamiltonian spaces\/} if there is an isotropy
submanifold $\Omega \subset X\times \overline{F}\times E$, such that
\begin{enumerate}
 \item $\Omega $ is a graph over both
$X\times_{G_{0}}\overline{F}$ and $X\times_{H_{0}}E$; and
 \item $(yx^{-1})\cdot f=y(x^{-1}(f))$ and $(x^{-1}z)\cdot
e=x^{-1}(z(e))$,
\noindent whenever either side is defined for any $x,y, z\in X, e\in E$
and $f\in F$, where by $x^{-1}(f)$ (or $z(e)$ resp.),
 we denote the unique
element in $E$ (or $F$ resp.) such that $(x, f, x^{-1}(f))\in \Omega $ 
(or $(z, z(e), e)\in \Omega $ resp.), and 
 $yx^{-1}$ (or $x^{-1}z$ resp.) denotes the
corresponding element $[y, x]$ (or $[x, z]$ resp.) in the groupoid
$G$ (or $H$ resp.) under the identification: 
 $G \cong (X\times_{H_{0}}X)/H $  (or $H\cong
 G\setminus (X\times_{G_{0}}X)$ resp.)
\end{enumerate}
\end{defn}

The following property follows immediately from the  definition above:

\begin{prop}
\begin{enumerate}
\item $x^{-1}(x(e))=e$ and $x(x^{-1}(f))=f$ for all composable $x\in X, e\in E$
and $f\in F$;
 \item for all composable $g\in G, x, y\in X, h\in H, f\in F$ and $e\in E$,
\be
&&(g\cdot x)^{-1} (f)  = x^{-1} (g^{-1} \cdot f),  \  (g\cdot y) (e)  = g \cdot
y(e);\\
&&(x\cdot h)^{-1} (f)  = h^{-1} \cdot (x^{-1}(f)),\   (y\cdot h) (e) 
= y(h \cdot e).
\ee
\end{enumerate}
\end{prop}

We are now ready to prove the main result of this section.

\begin{them}
\label{thm:bijection}
Suppose that $({G}\toto G_0 , \omega_G +\Omega_G )$
and $({H}\toto H_0 , \omega_H +\Omega_H )$ are Morita equivalent
quasi-symplectic groupoids with
  an equivalence bimodule $G_{0}\stackrel{\rho }{\leftarrow} X
\stackrel{\sigma}{\rightarrow} H_{0}$.  Then,
\begin{enumerate}
\item corresponding to
any Hamiltonian $G$-space $\phi:F \to G_0$, there is a unique (up to isomorphism)
Hamiltonian $H$-space $\psi:E \to H_0$
 such that $F$ and $E$ are a pair of related Hamiltonian spaces
 and vice versa.
\item let $\phi_i :F_i \to G_0$, $i=1, 2$, be Hamiltonian $G$-spaces
and $\psi_i: E_i\to H_0$, $i=1, 2$, their related  Hamiltonian
$H$-spaces. If $\phi_1$ and $\phi_2$ are clean, then
$\psi_1$ and $\psi_2$ are clean, and the classical intertwiner
spaces $F_1\times_G\overline{F_2}$ and
$E_1\times_H \overline{E_2}$ are  symplectically diffeomorphic.
\end{enumerate}
\end{them}
\begin{pf}
The proof is a simple modification of Theorem 4.2 in \cite{Xu:90}.

(1) Suppose that $\phi:F \to G_0$ is an Hamiltonian $G$-space.
Then  $G_{0}\stackrel{\phi }{\leftarrow} F
\stackrel{}{\rightarrow} \cdot$ is an Hamiltonian $G$-$\cdot$-bimodule.
Since  $H_{0}\stackrel{\sigma }{\leftarrow} \overline{X}
\stackrel{\rho}{\rightarrow} G_{0}$ is an Hamiltonian
$H$-$G$-bimodule, from Theorem \ref{thm:bimodule}  it follows that
$E:=\overline{X}\times_G F$ is an Hamiltonian
 $H$-$\cdot$-bimodule, i.e.,
an Hamiltonian $H$-space. Here
 $\psi : E\to H_{0}$ and the  $H$-action on $E$ are defined by
$$
\psi ([x,f])=\sigma (x)
$$
and
$$
h\cdot [x, f]=[x\cdot h^{-1}, f].
$$

Let $\Omega =\{(x, f, [x,f])|\forall (x, f)\in {X}\times_{G_{0}}F\}
\subset X\times \overline{F}\times E$. It is straightforward
to check that $\Omega$ is an isotropy submanifold,
and  is indeed a graph over both $X\times_{G_{0}}F$ and
$X\times_{H_{0}}E$. Hence $\phi:F \to G_0 $ and $\psi:E \to H_0 $ 
are a pair of related Hamiltonian spaces. Conversely,  
one easily sees that $F\cong X\times_H E$ by working backwards.

(2). Let $\Omega_i \subset X\times \overline{F_i}\times E_i$, $i=1, 2$,
be as in (1). Then $\Omega_1 \times \overline{\Omega_2} 
\subset X\times \overline{F_1}\times E_1
\times \overline{X}\times {F_2}\times \overline{E_2} $
 is an isotropy submanifold,
which is a graph over 
$X\times_{G_0} \overline{F_1} \times \overline{X}\times_{G_0} {F_2}$.
Given any $[(f_1 , f_2)]\in F_1 \times_{G} F_2$,
take any element $x \in X$ such that $\rho (x)=\phi_1 (f_1 )=\phi_2 (f_2 )$.
Let $e_1 \in E_1$,  and $e_2 \in E_2$ such that
$(x, f_1 , e_1 , x, f_2 , e_2 )\in \Omega_1 \times \Omega_2 $.
Then it is simple to see that $(e_1 , e_2 )\in E_1 \times_{H_0} E_2$
and $[e_1 , e_2 ]\in  E_1 \times_{H}E_2$ is independent of
the choice of $x$ and $(f_1 , f_2)$. Thus, we obtain a well-defined map:
$$\Phi: F_1\times_G\overline{F_2}
\to E_1\times_H \overline{E_2}, \ \ \ \ [f_1 , f_2 ]\to [e_1, e_2 ]. $$
It is simple to check that $\Phi$ is a bijection,
which is indeed  a symplectic diffeomorphism by using
the fact that $\Omega_1 \times \overline{\Omega_2}   $ is isotropic.
\end{pf}

\begin{cor}
\label{cor:reduction}
Assume that $({G}\toto G_0 , \omega_G +\Omega_G )$
and $({H}\toto H_0 , \omega_H +\Omega_H )$ are Morita equivalent
quasi-symplectic groupoids, and $\phi: F\to G_0$
and $\psi: E\to H_0$ are a pair of  related Hamiltonian
$G$- and $H$-spaces respectively.
Let $n\in H_{0}$ and  $m\in G_0$ be a pair
of related points. Then the
reduced spaces  $\phi^{-1}(m)/G_m^m$ and $\psi^{-1}(n)/H_n^n$
 are symplectically diffeomorphic.
\end{cor}

\begin{numrmk}
Corollary \ref{cor:reduction} indicates that the reduction of
Hamiltonian spaces of quasi-symplectic groupoids is of stack
natural.
\end{numrmk}

In fact, the  same  argument  in the proof of Theorem \ref{thm:bijection} leads to the following  more general result.


\begin{them}
Assume  that $f: G_{0}\stackrel{\rho }{\leftarrow} X
\stackrel{\sigma}{\rightarrow} H_{0}$
is a generalized  homomorphism of  quasi-symplectic groupoids
 from $({G}\toto G_0 , \omega_G +\Omega_G )$
to  $({H}\toto H_0 , \omega_H +\Omega_H )$. Then
\begin{enumerate}
\item if $\phi:E \to H_0$ is an Hamiltonian $H$-space,
and the maps $\phi$ and $\sigma$ are clean, then
 $\psi:F \to G_0$, where $F=X\times_H E$,  is an  Hamiltonian $G$-space,
called the pull-back Hamiltonian space and denoted by $f^*E$;
\item let $\phi_i :E_i \to H_0$, $i=1, 2$, be Hamiltonian $H$-spaces
and $\psi_i: F_i\to G_0$, $i=1, 2$, their pull-back
  Hamiltonian $G$-spaces. If $\phi_1$ and $\phi_2$ are clean, then
$\psi_1$ and $\psi_2$ are clean, and moreover
there exists a  natural symplectic immersion between
their  classical intertwiner spaces $F_1\times_G\overline{F_2}
\to E_1\times_H \overline{E_2}$.
\end{enumerate}
\end{them}

\subsection{Examples}


In this subsection, we will discuss various examples of
Morita equivalent quasi-symplectic groupoids and
derive some familiar corollaries as a consequence.  We start with a
general set-up.

Let  $(\gm\toto P, \omega +\Omega )$ be  a quasi-symplectic groupoid
and $\phi: Y\to P$ a surjective submersion. Consider
$(\gm [Y]\toto Y, \omega' +\Omega' )$, where $\omega'
=(\overline{B}, B, \omega )\in \Omega^2 (\gm [Y] ) $
and $\Omega' =\phi^* \Omega -dB$ (in applications,
normally $\phi^* \Omega=dB$ for some
$B\in \Omega^2 (Y)$,  so $\Omega' =0$).
 According to
Propositions \ref{prop:4.7} and \ref{prop:4.9},
 this is a  quasi-symplectic groupoid Morita equivalent
to $(\gm\toto P, \omega +\Omega )$.
Applying Theorem \ref{thm:bijection}, we obtain the following

\begin{prop}
\label{pro:5.1}
\begin{enumerate}
\item
There is a bijection between Hamiltonian $\gm$-spaces and
Hamiltonian $\gm [Y]$-spaces.

More precisely, if $(M\stackrel{J}{\to} P, \omega_M )$ is an
Hamiltonian $\gm$-space, then  $(N\stackrel{\tJ}{\to} Y, \omega_N )$
is  an Hamiltonian $\gm[Y]$-space, where $N$ is the fiber product
$Y\times_P M$, $\tJ: N\to Y$ is the projection
to the first component,  and $\omega_N =-\tJ^* B +p^* \omega_M$.  Here
$p: N\to M$ is the  projection to the second component.
 
Conversely, if $(N\stackrel{\tJ}{\to} Y, \omega_N )$ is  an Hamiltonian 
$\gm[Y]$-space, its corresponding 
$\gm$-space $(M\stackrel{J}{\to} P, \omega_M )$ is given as follows.
$M$ is the quotient space $N/\gm [Y]'$, where $ \gm [Y]'$ is  the
subgroupoid of $\gm [Y]$ consisting of all elements
$(y_1, y_2, u)$ with $y_1, y_2\in Y$, $\phi (y_1 )=\phi (y_2 )=u$,
$J:M\to P$ is given by
$J([n])=(\phi \smalcirc \tJ) (n)$,  and the two-form $\omega_M$ on $M$ is
defined by the equation: 
$$\pi^* \omega_M = \omega_N+ \tilde{J}^* B.$$ 
Here $\pi: N\to M$ denotes the natural  projection map.

\item If $(M\stackrel{J}{\to} P, \omega_M )$ and 
$(N\stackrel{\tJ}{\to} Y, \omega_N )$  are a pair of Hamiltonian
$\gm$- and $\gm [Y]$-spaces as above,
 and $\O\subset P$ and $\O_Y\subset Y$ are a pair
of related groupoid  orbits, then the reduced spaces
$J^{-1}(\O)/\gm$ and $\tJ^{-1}(\O_Y )/\gm [Y]$ are symplectic
diffeomorphic.
\end{enumerate}
\end{prop}
\begin{pf}
As in the proof of Propositions \ref{prop:4.7} and
 \ref{prop:4.9},
 the  Morita equivalence Hamiltonian
bimodule is given by  $P\stackrel{\rho }{\leftarrow} X
 \stackrel{\sigma }{\rightarrow} Y$, where $X= \gm \times_{t, P, \phi} Y$
and $\omega_X =( \omega, B)$.  The  left $\gm$- and right-
$\gm [Y]$-actions are given by 
Eqs. (\ref{eq:action1})-(\ref{eq:action2}) respectively.

Now we are ready to apply Theorem \ref{thm:bijection}.
If $J: M\to P$ is an Hamiltonian $\gm$-space, then
its corresponding Hamiltonian  $\gm [Y]$-space
is $N=\overline{X}\times_\gm M$, which is the quotient by $\gm$ of the space
$\{(r, y, m)|t(r)=\phi (y), \ J(m)=s (r)\}$.
It is simple to see that the latter is  diffeomorphic to the fiber product
$Y\times_P M$, and, under this diffeomorphism,
 the two-form on $\overline{X}\times_\gm M$
goes to $  -\tJ^* B +p^* \omega_M$.

Conversely, assume that $\tJ: N\to Y$ is  an  Hamiltonian $\gm [Y]$-space.
Then $M=X\times_{\gm [Y]}N\cong (X\times_{Y}N) /\gm [Y]$.
Now $X\times_{Y}N=\{(r, y, n)|t(r)=\phi (y), \ \tJ (n)=y\}$.
It is simple to see that, under the $\gm [Y]$-action,
any element in $X\times_{Y}N$ is equivalent to  $(u, y, n)$
 where $y=\tJ (n) $ and $u=\phi (\tJ (n))$.
Any two such elements $(u, y, n)$ and $(u', y', n')$ are equivalent
if and only if $n' =\gamma'  \cdot  n $ where $\gamma'\in  
\gm [Y]'$. As a result, $M$ can be identified with 
$N/\gm [Y]'$, and the two-form $(\omega, B, \omega_N)$ on $X\times_{Y}N$
goes to $\omega_N +\tilde{J}^* B$ under the identification
$$\{(u, y, n)|\forall n\in N, \ y=\tJ (n), \ u=\phi (\tJ (n))\} \longiso N.$$
 Therefore we have $\pi^* \omega_M = \omega_N+ \tilde{J}^* B$.

The rest of the claims follows easily from Theorem \ref{thm:bijection}.
\end{pf}

We  now  consider various special cases of the above proposition.

Let $G$ be a compact connected Lie group equipped with the
Bruhat-Poisson group
 structure \cite{LuW}, and $\frakg$ be its Lie algebra.
By $G^*$ we denote its simply-connected dual Poisson group.
It is known that there exists a  diffeomorphism \cite{A, AMW}:
$$E: \frakg^* \to G^*, $$
which is $G$-equivariant with respect to the coadjoint action 
on $\frakg^*$ and the left dressing action on $G^*$.
Let us recall the construction briefly. Here we follow
the presentation of \cite{AMW}.

Let $\kappa:\,\frakg^\complex\to\frakg^\complex$ be the
 Cartan involution given 
by  the complex conjugation,
 and let $\dagger:\,\frakg^\complex\to\frakg^\complex$ be the anti-involution
$ \xi^\dagger = -\kappa (\xi)$.
We also denote by $\dagger$ the induced anti-involution of $G^\complex$,
considered as a real group.  
Let $B^\sharp: \frakg^* \rightarrow \frakg$ be
the isomorphism induced by the Killing form $B$.
 For any $\mu \in \frakg^*$, the element $g=
\exp(i B^\sharp(\mu)) \in G^\complex $ admits a unique
decomposition $g = {l}{l}^\dagger$, for some ${l}
\in G^*$.  Then $E$ is defined by $E (\mu  )=l$.

Let $\beta \in \Omega^1 (\frakg^* )$ be the one-form \cite{AMW}
\begin{equation}\label{eq:beta}
\beta=\frac{1}{2i} \mathcal{H}\big( E^* B^\complex (\theta, \theta^\dagger) \big),
\end{equation}
where    $\theta\in \Omega^1 (G^*)\otimes\frakg^*$
 is the left-invariant Maurer-Cartan form, and  $\theta^\dagger$  its image
 under the map $\dagger:\,\frakg^*\subset\frakg^\complex \to\frakg^\complex$,
$\mathcal{H}:\,\Omega^\star(\frakg^*)\to \Omega^{\star-1}(\frakg^*)$
is the standard homotopy operator for the de Rham differential.
Let $B=d\beta\in \Omega^2 (\frakg^* )$.

The following proposition also follows from
Ginzburg-Weinstein theorem \cite{GW}.

\begin{prop}
The Lu-Weinstein symplectic groupoid
$G\times G^* \toto G^* $ is Morita equivalent to the
standard cotangent symplectic groupoid $T^*G\toto \frakg^*$.            
\end{prop}
\begin{pf}
Since $E: \frakg^* \to G^* $ is $G$-equivariant, the pull-back
groupoid $(G\times  G^* ) [\frakg^* ]$ is clearly isomorphic
to the transformation groupoid $G\times \frakg^*\toto \frakg^*$,
which is naturally  isomorphic to $T^*G\toto \frakg^*$.
Moreover, from  Lemma 2 (2) in \cite{A}
(or Proposition 3.1 in \cite{AMW}), it follows  that
$$E^* \omega' -\omega =\partial B.$$
Therefore, these two symplectic groupoids are Morita equivalent
since $dB=0$.
\end{pf}

As an application, we are lead to the following
 Alekseev-Ginzburg-Weinstein linearization theorem \cite{A}.

\begin{cor} 
Let $G$ be a  connected compact Lie group equipped with the
Bruhat-Poisson group structure. Then
\begin{enumerate}
\item   $(M , \omega_M )$ is
an Hamiltonian Poisson group
 $G$-space with  the momentum map $J: M\to G^* $
if and only if  $(M , \omega'_M )$
is a  usual  Hamiltonian $G$-space with the momentum map
$\tJ: M\to \frakg^*$, where
$$ J = E\smalcirc \tJ, \ \ \omega'_M=\omega_M -\tJ^*B. $$
\item If $\tilde{\O}$ is a coadjoint orbit in $\frakg^*$ and
$\O =E(\tilde{\O})$ is its corresponding dressing orbit
in $G^*$, then the reduced spaces $\tJ^{-1}(\tilde{\O})/G$ and
$J^{-1}( \O )/G$ are symplectically diffeomorphic.
\end{enumerate} 
\end{cor} 

Next we consider the AMM quasi-symplectic groupoid
$(G\times G\toto G, \omega +\Omega )$. Let 
 $\hol: L\frakg\lon G$ be the holonomy map, i.e., the time-$1$ 
 map of the differential equation:  
$$\hol_s (r)^{-1} \frac{\partial}{\partial s} \hol_s (r) =r, \ \hol_0 (r)=e. $$   Then we have $\hol^* \Omega  =d\mu$, where $\mu $ is the two-from
 on $L\frakg$ \cite{AMM}:
$$\mu=\half \int_0^1 \la \hol_s^* \btheta , \frac{\partial}{\partial s} \hol_s^* \btheta )ds, $$ 
where $\btheta \in \Omega^1 (G)\otimes \frakg$ is the right 
Maurer-Cartan form.

The pull-back groupoid of the AMM-groupoid under the 
holonomy map is isomorphic to the transformation groupoid
$LG\times L\frakg \toto { L\frakg}$, where $LG$ acts on
$L\frakg$ by the gauge transformation (\ref{eq:gauge}). 
 To see this, note that 
$$(G\times G ) [L\frakg ]\cong 
\{(r_1 (s), r_2 (s)  , g  )|r_1 (s) , r_2  (s)\in L\frakg , g
\in G \mbox{ such that } g^{-1} \hol (r_1 )g =\hol (r_2 )\}$$

Define 
\begin{equation}
\label{eq:GG}
\tau: \ (G\times G ) [L\frakg ]\to LG\times L\frakg , \ \ \ (r_1 (s),
r_2 (s) , g )
\to (r_1 (s), g(s)), 
\end{equation}
where $g(s) $ is defined by

\begin{equation}
Ad_{g(s)^{-1}} r_1 (s)- g(s)^{-1}\frac{dg(s)}{ds}=r_2 (s), \ \ g(0)=g.
\end{equation}
 
It is simple to see that $\tau $ is indeed a diffeomorphism, under which
the groupoid structure on $(G\times G ) [L\frakg ]$ becomes
the  transformation groupoid $LG \times L\frakg \toto L\frakg $.

\begin{prop}
\cite{BXZ}
The symplectic groupoid $(LG\times L\frakg \toto { L\frakg},
\omega_{LG \times L\frakg} )$ is Morita
equivalent to the AMM quasi-symplectic groupoid
$(G\times G\toto G, \omega+\Omega)$.
\end{prop}
\begin{pf}
From the above  discussion, we know that $LG\times L\frakg \toto { L\frakg}$
is the pull-back groupoid of $G\times G\toto G$ under the holonomy
map $\hol$. Denote by $f$ the groupoid homomorphism
from $LG \times L\frakg \toto L\frakg$ to $G\times G\toto G$, 
where on the space of morphisms and  the space of
objects, $f$ is given, respectively, by
$f (g(s), r(s))=(g(0), \hol (r) )$ and $f (r (s) )=\hol (r)$,
 $\forall g(s)\in
LG, \ r(s)\in L\frakg $. Then a simple computation yields that
$$\omega_{LG \times L\frakg}-f\upst (\omega +\Omega )=\delta \mu.$$ 
Thus the conclusion follows from Propositions \ref{prop:4.7}
and \ref{prop:4.9}  immediately.
\end{pf}

\begin{numrmk}
The above result was used in \cite{BXZ} to construct
an equivariant $S^1$-gerbe over the stack $G/G$.
\end{numrmk} 

An immediate consequence is the following equivalence theorem
of Alekseev--Malkin--Meinrenken \cite{AMM}.

\begin{cor}
\begin{enumerate}
\item There is a bijection between
 Hamiltonian $LG$-spaces  and quasi-Hamiltonian G-spaces. 

More precisely, if $(M\stackrel{J}{\to} G, \omega_M )$ is a
quasi-Hamiltonian $G$-space,
 then  $(N\stackrel{\tJ}{\to} L\frakg , \omega_N )$
is  an Hamiltonian $LG$-space, where $N$ is the fiber product
$ L\frakg\times_G M$, $\tJ: N\to  L\frakg$ is the projection map
to the first component,  and $\omega_N =-\tJ^* \mu +p^* \omega_M$.  Here
$p: N\to M$ is the  projection to the second component.
 
Conversely, if $(N\stackrel{\tJ}{\to}  L\frakg, \omega_N )$ is  an
 Hamiltonian $ LG$-space, its corresponding
quasi-Hamiltonian $G$-space $(M\stackrel{J}{\to} G, \omega_M )$ is given as follows.
$M$ is the quotient space $N/\Omega G$, where $ \Omega G$ is
the based loop group $\Omega G\subset LG$, $J: M\to G$ is given by
$J([n])=( \hol \smalcirc \tJ ) (n)$, and the two-form on $M$ is
defined  by 
$$\pi^* \omega_M = \omega_N+ \tJ^* \mu, $$ where $\pi: N\to M$ denotes
the projection.
\item Let $(M\stackrel{J}{\to} G, \omega_M )$ and
$(N\stackrel{\tJ}{\to} L\frakg, \omega_N )$  be as above.
Then the reduced spaces
$J^{-1}(e)/G$ and $\tJ^{-1}(0 )/LG $ are symplectically diffeomorphic.  
\end{enumerate}
\end{cor}
\begin{pf}
This essentially follows from Proposition \ref{pro:5.1}. Note that 
under the isomorphism (\ref{eq:GG}), the subgroupoid
$(G\times G) [L\frakg ]'$ of $(G\times G) [L\frakg ]$ corresponds to the
transformation  groupoid $L\frakg \times \Omega G\toto L\frakg$.
\end{pf}




\begin{numrmk}
For a quasi-Manin triple $(d, \frakg, \frakh )$,
Alekseev and Kosmann-Schwarzbach  introduced  a momentum
map theory with  target space $D/G$ \cite{AK-S}.
It would be interesting to investigate
what the corresponding quasi-symplectic groupoid is.
In particular, different choices of
 complements  $ \frakh$ should  give rise to
Morita equivalent quasi-symplectic groupoids.          
\end{numrmk}

\end{document}